\documentclass[12pt]{amsart}
\usepackage{amsfonts,latexsym,rawfonts,amsmath,amssymb,amsthm,times,a4wide}
\usepackage{geometry} 
\usepackage{amsmath}
\usepackage{amsthm}
\usepackage{amsfonts,amssymb}
\usepackage{fancyhdr} 

\usepackage{mathrsfs}  
\usepackage{graphicx}

\usepackage[all]{xy}

\usepackage{times}
\usepackage{cases}
\usepackage{color}

\usepackage[colorlinks,linkcolor=blue,anchorcolor=blue,citecolor=blue]{hyperref}

\usepackage{float}
\usepackage{subfig}
\usepackage{amsthm}

\newtheorem{question}{Question}
\newtheorem{conj}{Conjecture}

\newtheorem{thm}{Theorem}[section]
\newtheorem{lemma}[thm]{Lemma}
\newtheorem{prop}[thm]{Proposition}
\newtheorem{cor}[thm]{Corollary}

\newtheorem{rmk}{Remark}[section]

\theoremstyle{definition}
\newtheorem{definition}[thm]{Definition}

\newcommand{\IH}{{\mathbb H}}

\newcommand{\IK}{{\mathbb K}}

\newcommand{\IP}{{\mathbb P}}

\newcommand{\IR}{{\mathbb R}}

\newcommand{\CB}{{\mathcal B}}

\newcommand{\CD}{{\mathcal D}}

\newcommand{\CL}{{\mathcal L}}

\newcommand{\CR}{{\mathcal R}}

\newcommand{\CT}{{\mathcal T}}

\newcommand{\wt}{\widetilde}

\newcommand \di{\,\mathrm{d}}

\begin{document}

\title{Combinatorial Ricci Flow and Thurston's Triangulation Conjecture}
\author{Ke Feng, Huabin Ge}
\address{}
\curraddr{}
\email{}
\thanks{}
\keywords{Thurston's geometrical ideal triangulation, rigidity, hyperbolic tetrahedron, combinatorial Ricci flow}
\date{\today}
\dedicatory{}

\begin{abstract}
Thurston's triangulation conjecture asserts that every hyperbolic 3-manifold admits a geometric decomposition into ideal hyperbolic tetrahedra, a result proven only for certain special 3-manifolds. This paper presents combinatorial Ricci flow as a systematic and general approach to addressing Thurston's triangulation conjecture, showing that the flow converges if and only if the triangulation is geometric. First, we prove the rigidity of the most general hyperbolic polyhedral 3-manifolds constructed by isometrically gluing partially truncated and decorated hyperbolic tetrahedra, demonstrating that the metrics are uniquely determined by cone angles modulo isometry and decoration changes. Then, we demonstrate that combinatorial Ricci flow evolves polyhedral metrics toward complete hyperbolic structures with geometric decompositions when convergent. Conversely, the existence of a geometric triangulation guarantees flow convergence.
\end{abstract}

\maketitle


\section{Introduction}
Based on Thurston's geometrization program, vast majority of compact 3-manifolds are hyperbolic. For cusped hyperbolic 3-manifolds, Thurston's creative idea is to obtain them by  pasting isometrically ideal hyperbolic tetrahedra (see \cite{CDGW2016,Thurston2022,Weeks1980,Martelli} e.g.). Conversely, it is natural to ask if it is possible to decompose each cusp 3-manifold into ideal hyperbolic tetrahedra. This remains an open question, and is referred to as Thurston's geometric ideal triangulation conjecture in our article, since Thurston once proved the famous hyperbolic Dehn filling theorem on the premise of admitting that it was correct. To some extent, it is one of the most unsolved problems in the field of 3D geometry and topology since the Thurston's Geometrization Conjecture and the Virtual Haken Conjecture were solved. Similar to the cusp case, for a hyperbolic 3-manifold with both cusps and totally geodesic boundary, we conjecture it can be decomposed geometrically into partially truncated hyperideal tetrahedra. These conjectures are uniformly referred to as Thurston's triangulation conjecture for short.

The first breakthrough on Thurston's triangulation conjecture is due to Epstein-Penner \cite{Epstein-Penner}. They introduced canonical geometric ideal polyhedral decompositions on cusped hyperbolic manifolds with finite volume. Similarly, Kojima \cite{kojima} obtained truncated hyperideal polyhedral decompositions on compact hyperbolic 3-manifolds with totally geodesic boundary. For any hyperbolic 3-manifold $M$, although a geometric ideal triangulation had been proven to exist virtually by Luo-Schleimer-Tillmann \cite{Luo-Schleimer}, ``it is a difficult problem in general", as was pointed by Gu\'{e}ritaud and Schleimer in \cite{GS2010}, ``General results are known only when $M$ is restricted to belong to certain classes of manifolds: punctured-torus bundles, two-bridge link complements, certain arborescent link complements and related objects, or covers of any of these spaces" (see \cite{Aki}-\cite{Aki-SWY-2}, \cite{Gue-1,Gue-2,Ham-P,Jorgen,Lackenby,Nimer} for example).

Inspired by using Ricci flow and geometric analysis methods to solve the famous Poincar\'{e} conjecture, Chow-Luo \cite{Chow-Luo} and Luo \cite{Luo2005} initiated a program for studying the geometry and topology of three-dimensional manifolds using combinatorial Ricci flows. Their program is being further developed and improved in \cite{Feng2022}, where the authors proved that, for a compact 3-manifold with no sphere or torus boundary components, any topological ideal triangulation with edge valence no less than 10 is geometric. The main idea of this program is to paste hyperbolic tetrahedra to obtain a hyperbolic polyhedral manifold $M$, and a topological triangulation $\CT$ with cone singularity on each edge of $\CT$. A central observation from Thurston is that the measure at which the cone angle is equal to $2\pi$ happens to be non-singular, and thus hyperbolic. Hence $M$ admits a complete hyperbolic metric and a geometric triangulation isotopic to $\CT$. The combinatorial Ricci flow provides a perfect tool for deforming the cone metric on $M$ to a metric without singularity. The final metric is hyperbolic and carries a geometric decomposition. In dealing with this program, the following basic guiding issue need to be addressed:
\begin{question}
\label{question-introduction}
Is there exists a hyperbolic polyhedral metric on $(M,\CT)$ with cone angle $2\pi$? Is it unique, and how to find it?
\end{question}

The uniqueness part, also known as rigidity, is relatively easy to handle. The basic building blocks to obtain a hyperbolic polyhedral 3-manifold are \emph{partially truncated and decorated hyperbolic tetrahedra}. The vertices of such a tetrahedron are either ideal or hyper-ideal. And it has a decoration, i.e. a horosphere centered at each ideal vertex and is truncated at each hyper-ideal vertex. Suppose $(M,\CT)$ is a triangulated compact pseudo 3-manifold with a triangulation $\CT$ and the set of edges $E$. By definition, a \emph{decorated hyperbolic polyhedral metric} is obtained by replacing each tetrahedra in $\CT$ by a partially truncated and decorated hyperbolic tetrahedron and replacing the affine glusing homeomorphisms by isometries preserving the decoration. These glued space are generally hyperbolic 3-manifolds with cone singularities. The curvature of the metric assigned to each edge $e\in E$ is equal to $2\pi$ minus the cone angle at $e$ for interior edge and $\pi$ minus the cone angle for boundary edge. By the construction, these polyhedral metrics are determined by the lengths of the edges. The following rigidity result completely solves the uniqueness part of the above Question \ref{question-introduction}.
\begin{thm}[Rigidity]
\label{Rigidity}
Suppose $(M,\CT)$ is a triangulated compact pseudo 3-manifold, then a decorated hyperbolic polyhedral metric on $(M,\CT)$ is determined up to isometry and change of decorations by its curvature.
\end{thm}

Compared to previous work in \cite{Luo2018} and \cite{Feng-Ge-Liu}, which deals with the rigidity of a space glued solely by one type (ideal, hyper-ideal, or $1$-$3$ type respectively) hyperbolic tetrahedron, our rigidity holds true for a space glued by maybe five different type tetrahedra (see Section \ref{section-five-tetrahedra}).

Next we come to the existence part. For a compact 3-manifold with boundary equipped with an ideal triangulation, Luo \cite{Luo2005} introduced the following combinatorial Ricci flow
\begin{equation}
\label{def-Ricci-flow}
dl_i/dt=K_i
\end{equation}
to study the existence of hyperbolic metrics, where $l_i$ is the length of the $i$-th edge and $K_i$ is the combinatorial curvature of the cone metric $l=(l_1,\cdots,l_m)$ at the $i$-th edge. As the flow evolves, the polyhedral metric may degenerate, resulting in a solution that may only exist for a finite time. To overcome this difficulty, we introduce the \emph{extended combinatorial Ricci flow}
\begin{equation}
\label{exkl}
dl_i(t)/dt=\wt{K}_i(l(t)),
\end{equation}
where $l(t)=(l_1(t),\cdots,l_m(t))\in \IR^E$ is a (may be degenerate) decorated hyperbolic polyhedral metric, and $\wt{K}_i$ is the extended curvature, for details see Section \ref{section-metric-curvature} and Section \ref{section-Extend-RCF}.

\begin{thm}[Fundamental Theorem of the Combinatorial Ricci Flow]
\label{Fundamental-Thm}
Suppose $M$ is an oriented compact 3-manifold with boundary, no component of which is a 2-sphere, $\CT$ is an ideal triangulation of $M$. Then for any initial data $l_0\in\IR^E$, the solution $\{l(t)\}\subset \IR^E$ to the extended combinatorial Ricci flow \eqref{exkl} is unique, and exists for all time $t\geq0$. Moreover, there is a hyperbolic structure on $M$ so that $\CT$ is isotopic to a geometric triangulation if and only if $l(t)$ converges to a decorated hyperbolic polyhedral metric.
\end{thm}

Thus the combinatorial Ricci flow provides a suitable tool for finding the unique hyperbolic polyhedral metric on $(M,\CT)$ with curvature $0$, i.e. the complete hyperbolic metric on $M-\partial_t$ ($\partial_t$ is the toral boundary components of $M$) and geometric triangulation $\CT$ on $M$. Carefully comparing Question \ref{question-introduction}, we see Theorem \ref{Fundamental-Thm} provides the right answer to the question ``how to find it", and meanwhile, transforms the ``existence" part into a convergence problem of the combinatorial Ricci flow \eqref{exkl}. By Lackenby's work \cite{Lackenby2000,Lackenby-AGT}, we know
\begin{cor}
\label{thm-converge-imply-topo}
Suppose $M$ is an oriented compact 3-manifold with boundary, no component of which is a 2-sphere, $\CT$ is an ideal triangulation of $M$. If the solution $l(t)$ to the extended Ricci flow \eqref{exkl} converges, then $M$ is irreducible, atoroidal and not Seifert fibred. Moreover, $\CT$ supports an angle structure and is strongly $1$-efficient and therefore $1$-efficient and $0$-efficient.
\end{cor}
For the definition of $0$-efficient and (strongly) $1$-efficient triangulations of 3-manifolds, and their applications to finiteness theorems, decision problems, algorithms, computational complexity and Heegaard splittings, we refer to \cite{Garoufalidis,Jaco,Kang}. Begin with any initial ideal triangulation $\CT$ and any initial metric $l(0)$, we deform the hyperbolic polyhedral metric $l(t)$ according to Luo's combinatorial Ricci flow (\ref{def-Ricci-flow}). As we know, finite time degenerations may occur. Before encountering the degeneracy, we change the triangulation $\CT$ by Pachner moves, and the combinatorial Ricci flow continues to evolve under the new triangulation. Such a flow is called \emph{a combinatorial Ricci flow with surgery}, and it not only evolves the polyhedral metric, but also evolves the ideal triangulation. The following conjecture is hopefully true.
\begin{conj}
\label{conjecture-topo-imply-converge}
Let $M$ be an oriented compact 3-manifold with boundary, no component of which is a 2-sphere. Suppose $M$ is irreducible, atoroidal and not Seifert fibred. Then after a finite number of surgeries, the combinatorial Ricci flow will converge. Thus a hyperbolic structure, and meanwhile, a geometric ideal triangulation on $M$ are obtained.
\end{conj}
\begin{rmk}
Obviously, Conjecture \ref{conjecture-topo-imply-converge} shows that Thurston's geometric triangulation conjecture holds. Furthermore, the proof of our conjecture will give a new proof of Thurston's hyperbolization theorem.
\end{rmk}

To prove Theorem \ref{Rigidity}, we follow the spirit of \cite{Luo2018} and \cite{Feng-Ge-Liu}. The major differences between this paper and \cite{Luo2018} come from the fact that the $2$-$2$ type and $3$-$1$ type hyperbolic tetrahedron and the $4$-$0$ type (or $0$-$4$ type) tetrahedron have completely different geometry. Similar to the $1$-$3$ type case, but more prominent: compared with the $0$-$4$ type hyperbolic tetrahedron, the two types of hyperbolic tetrahedra both do not have good symmetry of dihedral angles, and compared with the $4$-$0$ type hyperbolic tetrahedron, they do not have mutual determinism between edge lengths and dihedral angles. Thus, the key point is establishing a new argument for the relationship between dihedral angles and decorated edges

To prove the Fundamental Theorem of CRF, we derive the long-time existence and the uniqueness of the extended \emph{CRF} (an abbreviation for the Combinatorial Ricci Flow) by similar proof in \cite{Feng2022-1}. The key ingredient for the proof of the last part of Theorem \ref{Fundamental-Thm}, is that the extended combinatorial Ricci flow is the negative gradient flow of a convex function on $\IR^E,$ related to the $H$-function which is defined in Section \ref{section-H-function}. By the existence of a zero-curvature decorated metric, one can show that the $H$-function is proper on $\IR^E/\sim$ and the zero-curvature decorated metric is the unique critical point of the $H$-function. By the standard ODE theory, we obtain the exponential convergence of CRF to the global minimizer for any initial data.

The paper is organized as follows. In section \ref{section-geometry}, we study the geometry of decorated hyperbolic tetrahedra, and give explicit relations between the decorated edge lengthes and dihedral angles, i.e., the Cosine Law in a decorated hyperbolic tetrahedron.
In section \ref{section-extension}, we meticulously study the degeneration behavior of $2$-$2$ type and $3$-$1$ type hyperbolic tetrahedra. We extend the value of dihedral angles to degenerated cases, and further extend the domain and value of the volume and co-volume functions while maintain their convexity.
In section \ref{section-volume}, based on the Schl\"{a}fli formula of the volume functions, we study properties of the co-volume functions for decorated hyperbolic tetrahedra. Further, we extend the definition of the co-volume functions, and show the extended co-volume functions are convex functions about the decorated edge lengths.
In section \ref{section-rigidity}, the core content is to prove the global rigidity of polyhedral hyperbolic 3-manifolds, i.e. Theorem ~\ref{Rigidity}. On the other hand,  we also introduce the $H$-function by the extended co-volume function, which will play an important role in the next section.
In section \ref{section-CRF}, we will study the long-time existence, uniqueness and convergence of the solution of the extended CRF, so as to prove Theorem \ref{Fundamental-Thm}.

~

\noindent
\textbf{Acknowledgements:}
The first author is supported by NSFSC, no. 2023NSFSC1286.
The second author is supported by NSFC, no. 12341102.

\section {The geometry of hyperbolic tetrahedra}
\label{section-geometry}
In this section, we study the geometry of decorated hyperbolic tetrahedra, and give explicit relations between the decorated edge lengthes and dihedral angles, i.e., the Cosine Law in these partially truncated and decorated hyperbolic tetrahedra.

\subsection{Five types of hyperbolic tetrahedra}
\label{section-five-tetrahedra}
Roughly speaking, for each integer $0\leq k\leq 4$, a $k$-$(4-k)$ type \emph{partially truncated and decorated hyperbolic tetrahedron} is a hyperbolic tetrahedron with $k$ hyper-ideal vertices and $4-k$ ideal vertices, which has a decoration, i.e. a horosphere centered at each ideal vertex and is truncated at each hyper-ideal vertex. To be precise, consider the Klein's projective model for $\IH^3$. In this model, $\IH^3$ is identified to the open unit ball $\IK^3$ in  $\IR^3\subset \IR\IP^3$. Geodesics of $\IH^3$ then correspond to the intersection of straight lines of $\IR^3$ with $\IK^3$, and the totally geodesic planes in $\IH^3$ are the intersection of linear planes with $\IK^3$. Let $\mathscr{P}\subset \IR^3$ be a compact convex Euclidean tetrahedron such that the vertex $v_i,(i \in \{1,2,\dots, k\})$ lies in $\IR^3\backslash \overline{\IK^3}$ and the other vertices $v_j, (j\in\{k+1,\dots, 4\})$ lie on the boundary of $\IK^3$. We require each edge between two vertices $v_i$ and $v_j$ meets $\IK^3$. Let $C_i$ be the cone with the apex $v_i$ tangent to $\partial\IK^3$ and $\pi_i$ be the half-space not containing $v_i$ such that $\partial \pi_i \cap \partial \IK^3= C_i \cap \partial \IK^3$. Then a \emph{$k$-$(4-k)$ type partially truncated hyperbolic tetrahedron} is given by $\sigma:= \mathscr{P}\cap \cap_{i=1}^k\pi_i$, and each $v_i$, $1\leq i\leq k$ is called a hyper-ideal vertex of $\sigma$, and the other vertices $v_j$, $k+1\leq j\leq 4$ are called ideal vertices of $\sigma$. By definition, a \emph{$k$-$(4-k)$ type partially truncated and decorated hyperbolic tetrahedron} is a pair of $(\sigma, \{H_{k+1}, \dots, H_{4}\})$, where $\sigma$ is a $k$-$(4-k)$ type partially truncated hyperbolic tetrahedron and $H_j$ is a horosphere centered at the ideal vertex $v_j$ for $j\in\{k+1,\dots, 4\}$. We call $\{H_{k+1}, \dots, H_{4}\}$ the decorations of $\sigma$. Throughout the article, for a $k$-$(4-k)$ type partially truncated and decorated hyperbolic tetrahedron, we always assume that $v_1, \dots, v_k$ are hyper-ideal vertices with truncations $\{\Delta_{1}, \dots, \Delta_{k}\}$, and $v_{k+1}, \dots, v_4$ are ideal vertices with decorations $\{H_{k+1}, \dots, H_{4}\}$ respectively.

To simplify the notation and terminology, a $k$-$(4-k)$ type partially truncated and decorated hyperbolic tetrahedron is often abbreviate as a \emph{$k$-$(4-k)$ type hyperbolic tetrahedron}, or a \emph{$k$-$(4-k)$ type tetrahedron}. 
There are five types of hyperbolic tetrahedra in total, i.e, $0$-$4$, $1$-$3$, $2$-$2$, $3$-$1$ and $4$-$0$ type hyperbolic tetrahedra, which are referred to ``\emph{hyperbolic tetrahedra}" uniformly throughout this paper. Using these notations and terminologies, a $0$-$4$-type hyperbolic tetrahedron is an ideal tetrahedron with a horosphere centered at each vertex (with no truncations). A $4$-$0$-type hyperbolic tetrahedron is exactly the hyper-ideal tetrahedron (with no decorations). It is a compact convex polyhedron in $\IH^3$ that is diffeomorphic to a truncated tetrahedron in $\IR^3$ with four right-angled hyperbolic hexagonal faces. It should be noted that these five types of hyperbolic tetrahedra are all called ``hyper-ideal tetrahedra" by Bao-Bonahon \cite{Bao2002}, and they include ideal tetrahedra as a special case.


For a hyperbolic tetrahedron with ideal vertices, the length of the edge connecting an ideal vertex to other vertices is infinity.
To characterize a meaningful cutting and pasting of tetrahedra along such edges, Penner~\cite{Penner1987} introduced the concept of decorated edge lengthes. Taking a decorated $k$-$(4-k)$-type hyperbolic tetrahedron $(\sigma, \{H_{k+1}, \dots, H_{4}\})$ for example. For the edge $e_{ij}$ connecting two ideal vertices $v_i$ and $v_j$ with decorations $H_i$ and $H_j$ respectively, the signed edge length $l_{ij}$ is defined as follows.
The absolute value $\vert l_{ij} \vert $ is the distance between $ H_i \cap e_{ij}$ and $H_j \cap e_{ij}$,
so that $l_{ij} > 0$ if $ H_i $ and $ H_j$ are disjoint and $l_{ij} \leq 0$ if $H_i \cap H_j \neq \emptyset$.
Similarly, for the edge $e_{ij}$ connecting an ideal vertex $v_i$ with decoration $H_i$ to a hyper-ideal vertex $v_j$ (or the truncated triangle $\Delta_j$ without confusions), define the signed edge length $l_{ij}$ as the algebraic distance between $H_i \cap e_{ij}$ and $\Delta_j$, it is negative if and only if $\Delta_j$ is contained in $H_i$. For the edge $e_{ij}$ connecting two hyper-ideal vertices $v_i$ and $v_j$, let $l_{ij}$ be the distance between its truncated arc, i.e. the geodesic arc between $\Delta_i$ and $\Delta_j$. In this case, the edge length $l_{ij}$ is always positive. The dihedral angle at $e_{ij}$ is denoted by $\alpha_{ij}$. Obviously, if one only change the decoration at any ideal vertex, all the dihedral angles will not change however. As a convention, we assume $l_{ij}=l_{ji}$ and $\alpha_{ij}= \alpha_{ji}$. Moreover, as long as the character $e_{ij}$, $l_{ij}$ or $\alpha_{ij}$ appears, it is always require that $i\neq j$ and $i,j\in\{1,2,3,4\}$.

\begin{figure}[htbp]
	\centering
	\subfloat[$l_{ij}>0$]{\includegraphics[width=.28\columnwidth]{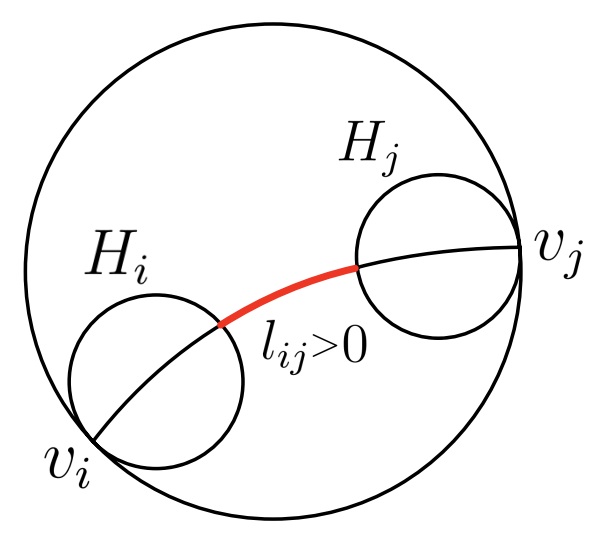}}\hspace{45pt}
	\subfloat[$l_{ij}<0$]{\includegraphics[width=.28\columnwidth]{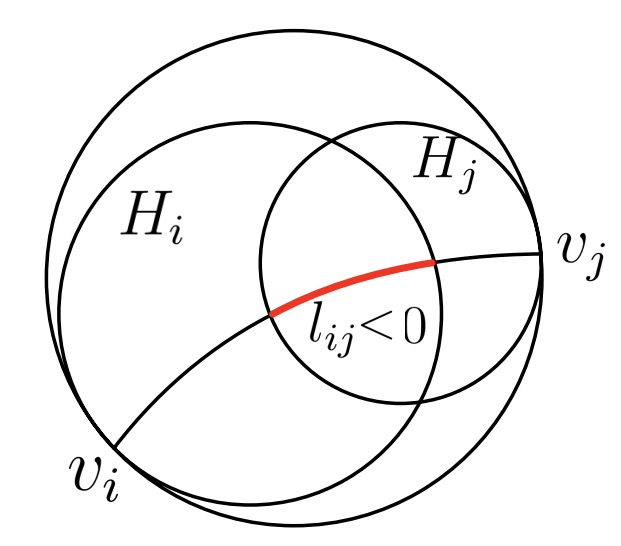}}
	\caption{decorated edge lengthes in the Poincar\'{e} disk model }
\end{figure}

Given a decorated $k$-$(4-k)$-hyperbolic tetrahedron $(\sigma, \{H_{k+1}, \dots, H_{4}\})$, $0\leq k\leq 4$, we fix the six edges in the order $(e_{12}, \dots, e_{34})=(e_{12}, e_{13}, e_{14}, e_{23}, e_{24}, e_{34})$ and fix the dihedral angles $(\alpha_{12}, \dots, \alpha_{34})$ and the decorated edge lengths $(l_{12}, \dots, l_{34})$ in the same order. The moduli space of all dihedral angles of decorated $k$-$(4-k)$-type hyperbolic tetrahedra is denoted by
\begin{equation}
\CB_{k|4-k} = \bigl\{(\alpha_{12}, \dots, \alpha_{34}) ~\text{is the dihedral angle vector of some $k$-$(4-k)$ type }~\sigma \bigr\}.
\end{equation}
Similarly, the moduli space of all decorated edge lengths of decorated $k$-$(4-k)$-type hyperbolic tetrahedra is denoted by
\begin{equation}
\label{equ-huati-L-k-4-k}
\CL_{k|4-k}=\bigl\{(l_{12}, \dots, l_{34})~\text{is the decorated edge length vector of some $k$-$(4-k)$ type }~\sigma \bigr\}.
\end{equation}
Preliminary, it's easy to see for each $0\leq k\leq 4$, the dihedral angle space $\CB_{k|4-k}\subset {\IR}^6_{>0}$ and
$$\CL_{k|4-k}\subset {\IR}^6.$$
It had been demonstrated \cite{Bao2002, Fujii1990, Leibon2002, MR1277810, Schlenker2002, Vinberg} that for a decorated hyperbolic polyhedron, all six dihedral angles can uniquely determine its shape, and all edge lengths can also uniquely determine its shape up to change of decorations. Although the moduli space $\CL_{k|4-k}$ is not easy to characterize directly, the moduli space $\CB_{k|4-k}$ has a very simple form of characterization. In fact, one have (for example see \cite{Bao2002} or Section 3 in \cite{Schlenker2002})
\begin{lemma}[]\label{prop 2.1}
Let $(\sigma, \{H_{k+1}, \dots, H_{4}\})$, $0\leq k\leq 4$ be a decorated $k$-$(4-k)$ type hyperbolic tetrahedron with hyper-ideal vertices $1, \dots, k$ and ideal vertices $k+1,\dots, 4$.
\begin{enumerate}
  \item[(a)] The dihedral angle vector $(\alpha_{12},\ldots, \alpha_{34})$ satisfies $\alpha_{ij}>0$, $\sum\limits_{j\neq i} \alpha_{ij}<\pi$ for each hyper-ideal vertex $i\in\{1,\dots,k\}$ and $\sum\limits_{j\neq i}\alpha_{ij}=\pi$ for each ideal vertex $i\in\{k+1, \dots, 4\}$.
  \item[(b)] The isometry class of $\sigma$ is uniquely determined by its six dihedral angles $(\alpha_{12},\ldots, \alpha_{34})$.
  \item[(c)] Given $(\alpha_{12},\ldots, \alpha_{34})\in{\IR}^6$ so that $\alpha_{ij}>0$, $\sum\limits_{j\neq i} \alpha_{ij}<\pi$ for each hyper-ideal vertex $i\in\{1,\dots,k\}$ and $\sum\limits_{j\neq i}\alpha_{ij}=\pi$ for each ideal vertex $i\in\{k+1, \dots, 4\}$, there exists a $k$-$(4-k)$ type hyperbolic tetrahedron having $\alpha_{ij}$ as its dihedral angles at the edge $e_{ij}$.
  \item[(d)] The isometry class of $\sigma$ is determined by its decorated edge length vector $(l_{12},\ldots,l_{34})$ up to change of decorations.
\end{enumerate}
\end{lemma}

Hence in $\sigma$, each dihedral angle is determined by the six decorated edge lengths, and the map $\alpha=\alpha(l): \CL_{k|4-k}\rightarrow\CB_{k|4-k}$, ~ $(l_{12},\ldots,l_{34})\mapsto(\alpha_{12},\ldots, \alpha_{34})$ is surjective and injective up to change of decorations. Consequently, there is a diffeomorphism between
$\CL_{k|4-k}$ (mod decorations) and $\CB_{k|4-k}$. In the following sections we will give explicit relations between the decorated edge lengthes and dihedral angles in $\sigma$, i.e., Penner's Cosine Law \cite{Penner1987} in $2$-$2$ type and $3$-$1$ type decorated hyperbolic tetrahedron respectively, and readers could get details for other cases (i.e. $1$-$3$ type, $0$-$4$ type and $4$-$0$ type tetrahedra respectively) in \cite{Feng-Ge-Liu} and \cite{Luo2018}.

\subsection {Geometry of $2$-$2$ type hyperbolic tetrahedra}
In this subsection, we study the geometry in $2$-$2$ type hyperbolic tetrahedra $(\sigma, \{H_3, H_4\})$. Here we give some notation. For one $2$-$2$ type hyperbolic tetrahedron, it has two ideal vertices $v_3, v_4$, two hyper-ideal vertices $v_1, v_2$ with truncated triangles $\Delta_1, \Delta_2$ respectively. Denote $P_{ijk} (i,j,k\in\{1,2,3,4\})$ by the hyperbolic polygonal face with edges (recall a vertex edge is not an ``edge") $e_{ij}$, $e_{ik}$ and $e_{jk}$. The vertex edge $\Delta_i \cap P_{ijk}$ is denoted
by $x^i_{jk}$ and the dihedral angle at $e_{ij} $ is denoted by $\alpha_{ij}$.
We also always assume $\alpha_{ij}= \alpha_{ji}$. Note each dihedral angle between every hyperbolic polygonal face and the vertex triangle is always $\pi/2$.

Thus, the space of isometry class of $2$-$2$  type tetrahedron parametrized by
dihedral angles is the open convex polytope in ${\IR}^6_{>0}$,
\[
    \CB_{2|2} =\bigl\{ (\alpha_{12},\ldots,\alpha_{34}) \in {\IR}^6_{>0} \mid \sum_{j\neq 1} \alpha_{1j} < \pi, \sum_{j\neq 2} \alpha_{2j} < \pi, \sum_{j\neq i} \alpha_{ij} = \pi \; \text{for}\; i\in \{3,4\}   \bigr\} .
\]


Now, for a decorated $2$-$2$ type tetrahedron with decorated edge length vector $(l_{12}, \dots, l_{34})\in {\IR}_{> 0} \times {\IR}^5 $,
let $\alpha_{ij}$ be the dihedral angle at edge $e_{ij}$, let $\theta_{ij}^k$ be the length
of $x^k_{ij}$ which is vertex edge or the intersection part of horospheres $H_k$ and the face $H_{kij}$. Similarly, based on the formulas in \cite{Guo2009}), we have the following equations,  for $\{i,j\}=\{1,2\}$ and $\{k,h\}=\{3,4\}$:

\begin{equation}\label{2-2edge-1}
\cosh \theta_{2k}^1=\frac{e^{l_{2k}}+e^{l_{1k}}\cosh l_{12}}{e^{l_{1k}}\sinh l_{12}}
\end{equation}

\begin{equation}
\cosh \theta_{1k}^2=\frac{e^{l_{1k}}+e^{l_{2k}}\cosh l_{12}}{e^{l_{2k}}\sinh l_{12}}
\end{equation}

\begin{equation}
\cosh \theta_{34}^i=2e^{l_{34}-l_{i3}-l_{i4}}+1
\end{equation}

\begin{equation}
(\theta_{12}^k)^2=\frac{2[\cosh l_{12}+\cosh(l_{1k}-l_{2k})]}{e^{l_{1k}+l_{2k}}}
\end{equation}

\begin{equation}
(\theta_{i3}^4)^2=\frac{e^{l_{i3}} + e^{l_{34}-l_{i4}}}{e^{l_{i4}+l_{34}}}
\end{equation}

\begin{equation}\label{2-2edge-6}
(\theta_{i4}^3)^2=\frac{e^{l_{i4}} + e^{l_{34}-l_{i3}}}{e^{l_{i3}+l_{34}}}
\end{equation}

\begin{eqnarray}
  \cos \alpha_{ik} &=& \frac{\cosh \theta_{jk}^i\; \cosh \theta_{kh}^i- \cosh \theta_{jh}^i}{\sinh \theta_{jk}^i \; \sinh\theta_{kh}^i} \notag \\[3pt]
  &=& \frac{(\theta_{ih}^k)^2+(\theta_{ij}^k)^2 - (\theta_{jh}^k)^2}{2\theta_{ih}^k \theta_{ij}^k},
\end{eqnarray}
and

\begin{eqnarray}
  \cos \alpha_{12} &=& \frac{\cosh \theta_{23}^1\; \cosh \theta_{24}^1 - \cosh \theta_{34}^1}{\sinh \theta_{23}^1 \; \sinh\theta_{24}^1} \notag \\[3pt]
  &=&  \frac{\cosh \theta_{13}^2\; \cosh \theta_{14}^2 - \cosh \theta_{34}^2}{\sinh \theta_{13}^2 \; \sinh\theta_{14}^2} ,
\end{eqnarray}

and
\begin{eqnarray}\label{2-2angle-1}
  \cos \alpha_{34} &=& \frac{(\theta_{14}^3)^2+(\theta_{24}^3)^2 - (\theta_{12}^3)^2}{2\theta_{14}^3 \theta_{24}^3} \notag\\[3pt]
  &=&  \frac{(\theta_{13}^4)^2+(\theta_{23}^4)^2 - (\theta_{12}^4)^2}{2\theta_{13}^4 \theta_{23}^4}.
\end{eqnarray}

\subsection {Geometry of $3$-$1$ type hyperbolic tetrahedra}
Next, we study the geometry in $3$-$1$ type hyperbolic tetrahedra $(\sigma, \{H_4\})$. For one $3$-$1$ type hyperbolic tetrahedron, it has one ideal vertex $v_4$, three hyper-ideal vertices $v_1, v_2,v_3$ with truncated triangles $\Delta_1, \Delta_2, \Delta_3$ respectively. Denote $P_{ijk} (i,j,k\in\{1,2,3,4\})$ by the hyperbolic polygonal face with edges (recall a vertex edge is not an ``edge") $e_{ij}$, $e_{ik}$ and $e_{jk}$. The vertex edge $\Delta_i \cap P_{ijk}$ is denoted
by $x^i_{jk}$ and the dihedral angle at $e_{ij} $ is denoted by $\alpha_{ij}$.
We also always assume $\alpha_{ij}= \alpha_{ji}$. Note each dihedral angle between every hyperbolic polygonal face and the vertex triangle is always $\pi/2$.

Thus, the space of isometry class of $M_{3|1}$-type tetrahedron parametrized by
dihedral angles is the open convex polytope in ${\IR}^6_{>0}$,
\[
    \CB_{3|1} =\bigl\{ (\alpha_{12},\ldots,\alpha_{34}) \in {\IR}^6_{>0} \mid \sum_{j\neq 4} \alpha_{j4} = \pi,  \sum_{j\neq i} \alpha_{ij} < \pi \; \text{for}\; i\in \{1,2,3\}   \bigr\} .
\]


Now, we also give some computational formulas between dihedral angles and decorated edge lengthes based on the formulas in \cite{Guo2009}.
For a decorated $3$-$1$ type tetrahedron with decorated edge length vector $(l_{12}, \dots, l_{34})\in {\IR}^2_{>0} \times {\IR} \times {\IR}_{>0} \times {\IR}^2$,
let $\alpha_{ij}$ be the dihedral angle at edge $e_{ij}$, let $\theta_{ij}^k$ be the length
of $x^k_{ij}$ which is vertex edge or the intersection part of horospheres $H_4$ and the face $H_{ij4}$, for $\{i,j,k\}=\{1,2,3\}$, we have the following equations:
\begin{equation}\label{3-1edge-1}
\cosh \theta_{j4}^i=\frac{e^{l_{j4}}+e^{l_{i4}}\cosh l_{ij}}{e^{l_{i4}}\sinh l_{ij}},
\end{equation}

\begin{equation}
(\theta_{ij}^4)^2=\frac{2[\cosh l_{ij}+\cosh(l_{i4}-l_{j4})]}{e^{l_{i4}+l_{j4}}},
\end{equation}

\begin{equation}\label{3-1edge-3}
\cosh \theta_{jk}^i=\frac{\cosh l_{jk}+\cosh l_{ij}\cosh l_{ik}}{\sinh l_{ij} \sinh l_{ik}}
\end{equation}
\begin{eqnarray}
  \cos \alpha_{i4} &=& \frac{\cosh \theta_{j4}^i\; \cosh \theta_{k4}^i - \cosh \theta_{jk}^i}{\sinh \theta_{j4}^i \; \sinh\theta_{k4}^i} \notag\\[3pt]
  &=&  \frac{(\theta_{ij}^4)^2+(\theta_{ik}^4)^2 - (\theta_{jk}^4)^2}{2\theta_{ij}^4 \theta_{ik}^4},
\end{eqnarray}

and
\begin{eqnarray}\label{3-1angle-1}
  \cos \alpha_{ij} &=& \frac{\cosh \theta_{j4}^i\; \cosh \theta_{jk}^i - \cosh \theta_{k4}^i}{\sinh \theta_{j4}^i \; \sinh\theta_{jk}^i} \notag \\[3pt]
  &=& \frac{\cosh \theta_{i4}^j\; \cosh \theta_{ik}^j - \cosh \theta_{k4}^j}{\sinh \theta_{i4}^j \; \sinh\theta_{ik}^j}.
\end{eqnarray}

\subsection{No lower dimensional degeneracy}
Similar to the argument in subsection 2.4 in \cite{Feng-Ge-Liu}, we also could get the same results on the degeneracy of $2$-$2$ type and $3$-$1$ type hyperbolic tetrahedra, i.e. there are no one-dimensional or two-dimensional degeneracy for $2$-$2$ type and $3$-$1$ type hyperbolic tetrahedra. More specifically, we have the following theorem.

\begin{thm}
\label{thm-non-degenerate}
Let $\CR$ be ${\IR}_{>0} \times {\IR}^5$ (or ${\IR}^2_{>0} \times {\IR} \times {\IR}_{>0} \times {\IR}^2$, respectively) when considering $2$-$2$ type (or $3$-$1$ type, respectively) tetrahedra. There is no one or two dimensional degeneracy for each $l\in \CR$, that is, for each $l\in \CR$, the two dimensional geometric configurations, i.e.
\begin{enumerate}
  \item in case of $2$-$2$ type: the two quadrilateral corresponding to $v_1v_3v_4$ and $v_2v_3v_4$, and the two pentagons corresponding to $v_1v_2v_3$ and $v_1v_2v_4$ (with corresponding decorations and truncations)
  \item in case of $3$-$1$ type: the three pentagons corresponding to $v_iv_jv_4$ ($\{i, j\}\subset\{1, 2, 3,\}$), and the triangle corresponding to $v_1v_2v_3$ (with corresponding decorations and truncations)
\end{enumerate}
always exist, and hence those $\{\theta^r_{pq}\}$ ($\{p, q, r\}\subset\{1, 2, 3, 4\}$) are always meaningful and can be obtained by the formulas (\ref{2-2edge-1})-(\ref{2-2edge-6}) or (\ref{3-1edge-1})-(\ref{3-1edge-3}) respectively.
\end{thm}

\section{The extension of the dihedral angles of hyperbolic tetrahedra}
\label{section-extension}
\subsection{$\phi$-functions and the decorated edge length space $\CR$}
As the necessary preparation, for $2$-$2$ and $3$-$1$ type tetrahedra, we define the $\phi$ functions respectively in this section, which are extensions of the cosine value of the dihedral angles.

\begin{definition}
\label{def-phi-22}
Consider the $2$-$2$ type case, let $l=(l_{12},l_{13},l_{14},l_{23},l_{24},l_{34})\in \CR=\IR_{>0}\times\IR^5$, and recall $l_{ij}=l_{ji}$, we define the functions $\phi_{pq}(l)$, $l\in\IR_{>0}\times\IR^5$ for $p, q\in\{1, 2,3,4\}$, $p\neq q$ as follows: set $\phi_{pq}(l)=\phi_{qp}(l)$ and
\begin{equation}\label{def-phi-12-22}
  \phi_{12}(l)=
  \frac{(e^{l_{13}}+e^{l_{23}}\cosh l_{12})(e^{l_{14}}+e^{l_{24}}\cosh l_{12})-(2e^{l_{34}}+e^{l_{23}+l_{24}})\sinh ^2 l_{12}}{[(e^{2l_{23}}+2e^{l_{23}+l_{13}}\cosh l_{12}+e^{2l_{13}})(e^{2l_{14}}+2e^{l_{24}+l_{14}}\cosh l_{12}+e^{2l_{24}})]^{1/2}},
\end{equation}
\begin{equation}\label{def-phi-34-22}
\phi_{34}(l)=\frac{e^{l_{24}+l_{13}}+e^{l_{14}+l_{23}}-2e^{l_{34}}\cosh l_{12}}{2[(e^{l_{24}+l_{23}}+e^{l_{34}})(e^{l_{14}+l_{13}}+e^{l_{34}})]^{1/2}},
\end{equation}
and for $\{i,j\}=\{1,2\}$, $\{k,h\}=\{3,4\}$, further define 
\begin{equation}
\phi_{ik}(l)=\frac{(2e^{l_{kh}-l_{ik}-l_{ih}}+1)(e^{l_{jk}}+e^{l_{ik}}\cosh l_{ij})-(e^{l_{jh}+l_{ik}-l_{ih}}+e^{l_{ik}}\cosh l_{ij})}{2[(e^{2(l_{kh}-l_{ik}-l_{ih})}+e^{l_{kh}-l_{ik}-l_{ih}})(e^{2l_{jk}}+2e^{l_{ik}+l_{jk}}\cosh l_{ij}+e^{2l_{ik}})]^{1/2}}.
\end{equation}
\end{definition}

\begin{lemma}\label{angleedge22}
The right hand side of (\ref{def-phi-12-22}) ((\ref{def-phi-34-22}), respectively) is symmetry about the lower index $1$ and $2$ ($3$ and $4$, respectively). Hence the setting $\phi_{12}(l)=\phi_{21}(l)$ and $\phi_{34}(l)=\phi_{43}(l)$ is well defined in Definition \ref{def-phi-22}. And if $l \in \CL_{2|2}$, that is, $l$ is the decorated edge length vector of a $2$-$2$ type tetrahedra $(\sigma, \{H_3, H_4\})$, then $\phi_{pq}(l)=\cos\alpha_{pq}$ for each pair $\{p,q\}\subset\{1, 2,3,4\}$.
\end{lemma}

\begin{proof}
Using (\ref{2-2edge-1})-(\ref{2-2angle-1}), $\alpha_{pq}(l)$ is determined by $l$ and hence a function of $l$. Moreover, by direct calculations, the symmetries and the equalities $\phi_{pq}(l)=\cos({\alpha_{pq}}(l))$ could be checked.
\end{proof}


\begin{lemma}\label{L-zero-22}
In the $2$-$2$ type case, the function $\phi_{ij}$ extends continuously to $\overline{\CR}={\IR}_{\ge 0}\times {\IR}^5 $. Furthermore,
  $\phi_{12}(l) = 1$ when $ l_{12} = 0 $.
\end{lemma}
\begin{proof}
   From Definition \ref{def-phi-22}, the denominator of $\phi_{ij}$ never equal to 0. Therefore, the function $\phi_{ij}$ can
   continuously extends ${\IR}_{\ge 0}\times {\IR}^5 $.
   Furthermore, a direct calculation show if $ l_{12} = 0 $, then $\phi_{12}(l) = \frac{(e^{l_{13}}+e^{l_{23}})(e^{l_{14}}+e^{l_{24}})}{(e^{l_{13}}+e^{l_{23}})(e^{l_{14}}+e^{l_{24}})} = 1$.
\end{proof}

\begin{definition}\label{def-phi-31}
Consider the $3$-$1$ type case, let $l=(l_{12}, \dots, l_{34})\in\CR={\IR}^2_{>0}\times {\IR}\times{\IR}_{>0}\times{\IR}^2$.
Define the functions $\phi_{pq}(l)$, $l\in{\IR}^2_{>0}\times {\IR}\times{\IR}_{>0}\times{\IR}^2$ as follows: set $\phi_{pq}(l)=\phi_{qp}(l)$ and
\begin{equation}
  \phi_{i4}(l)=
  \frac{e^{l_{j4}}e^{l_{k4}}+e^{l_{i4}}(e^{l_{k4}}c_{ij} + e^{l_{j4}}c_{ik}-e^{l_{i4}}c_{jk})}{[(e^{2l_{j4}}+2e^{l_{j4}+l_{i4}}c_{ij}+e^{2l_{i4}})(e^{2l_{k4}}+2e^{l_{k4}+l_{i4}}c_{ik}+e^{2l_{i4}})]^{1/2}};
\end{equation}
\begin{equation}
\label{def-phi-ij-31}
\phi_{ij}(l)=\frac{(e^{l_{j4}}+e^{l_{i4}} c_{ij})(c_{jk} + c_{ij} c_{ik})-(e^{l_{k4}}+e^{l_{i4}}c_{ik})s^2 _{ij}}{[(c^2 _{jk} + c^2_{ij} + c^2_{ik} + 2c_{jk}c_{ij}c_{ik}-1)(e^{2l_{j4}}+2e^{l_{j4}+l_{i4}}c_{ij}+e^{2l_{i4}})]^{1/2}}
\end{equation}
for $\{i,j,k\}=\{1,2,3\}$, where $c_{ij} = \cosh l_{ij}$, $s_{ij} = \sinh l_{ij}$.
\end{definition}

\begin{lemma}
\label{angleedge31}
The right hand side of (\ref{def-phi-ij-31}) is symmetry about the index $i$ and $j$, hence the setting $\phi_{ij}(l)=\phi_{ji}(l)$ ($\{i,j\}\subset\{1,2,3\}$) is well defined in Definition \ref{def-phi-31}. And if $l\in \CL_{3|1}$, that is, $l$ is the decorated edge length vector of a $3$-$1$ type tetrahedra $(\sigma, \{H_4\})$,
then $\phi_{pq}(l)=\cos\alpha_{pq}$ for each pair $\{p,q\}\subset\{1, 2,3,4\}$.
\end{lemma}

\begin{proof}
Using (\ref{3-1edge-1})-(\ref{3-1angle-1}), $\alpha_{pq}(l)$ is determined by $l$ and hence a function of $l$. Moreover, by direct calculations, the symmetries and the equalities $\phi_{pq}(l)=\cos({\alpha_{pq}}(l))$ could be checked.
\end{proof}

\begin{lemma}\label{L-zero-31}
In the $2$-$2$ type case, each function $\phi_{pq}$ ($1\leq p<q\leq4$) extends continuously to $\overline{\CR}={\IR}^2_{\ge 0} \times {\IR} \times {\IR}_{\ge 0} \times {\IR}^2$. Furthermore, for $\{i,j\} \subset\{2,3,4\}$, $\phi_{ij}(l) = 1$ when $ l_{ij} = 0 $.
\end{lemma}
\begin{proof}
   From Definition \ref{def-phi-31}, the denominators of all $\phi_{ij}$ never equal to 0. Therefore, the function $\phi_{ij}$ can continuously extends ${\IR}^2_{\ge 0} \times {\IR} \times {\IR}_{\ge 0} \times {\IR}^2$.
   Furthermore, a direct calculation show if $ l_{ij} = 0, (i,j \in {2,3,4})$, then $\phi_{ij}(l) = \frac{(e^{l_{i4}}+e^{l_{j4}})(\cosh l_{jk} + \cosh l_{ik})}{(e^{l_{i4}}+e^{l_{j4}})(\cosh l_{jk} + \cosh l_{ik})} = 1$
\end{proof}

\begin{rmk}\label{unchange31}
Given a $2$-$2$ type tetrahedron $(\sigma, \{H_1, H_2\})$, or a $3$-$1$ type tetrahedron $(\sigma, \{H_4\})$, we see from geometric intuition that as the radius of the decorated horosphere $H_i$ changes, although some decorated edge lengthes may change, the six dihedral angles do not change accordingly. This fact can also be easily checked from the definition of $\phi_{pq}$, Lemma \ref{angleedge22} and Lemma \ref{angleedge31}.
\end{rmk}

Recall $\CL_{k|4-k}$ is the space of all decorated edge lengths of a $k$-$(4-k)$-type tetrahedra. In the following, we use a consistent notation $\CL$ to represent these spaces, and no longer distinguish between $\CL$ and those $\CL_{k|4-k}$ at different cases.

\begin{prop}\label{L-space}
The space of all $2$-$2$ type (or $3$-$1$ type, respectively) hyperbolic tetrahedron parameterized by the edge lengths is
\[
  \CL = \{(l_{12},\ldots,l_{34}) \in {\CR} \mid \phi_{ij}(l)\in(-1,1) \; \text{for all } \; \{i,j\} \subset  \{1,2,3,4\}\},
\]
where $\CR$ equals to ${\IR}_{>0} \times {\IR}^5$ (or ${\IR}^2_{>0} \times {\IR} \times {\IR}_{>0} \times {\IR}^2$, respectively).
\end{prop}

\begin{proof}
We have proved this fact for $1$-$3$ type hyperbolic tetrahedra in \cite{Feng-Ge-Liu}. Here we only prove this result for $2$-$2$ type tetrahedra. The proof for $3$-$1$ type tetrahedra is similar.

If $l$ is the edge length vector of a $2$-$2$ type tetrahedron, then by definition, each dihedral angle $\alpha_{ij}(l) \in (0, \pi)$. Using the cosine law to the vertex triangles $\triangle_1 $, $\triangle_2 $ and polygons $H_{123}, \, H_{124}, \, H_{134}$, $H_{234}$, we have $\phi_{ij}(l) = \cos(\alpha_{ij}) \in (-1,1)$ for all
$\{i,j\} \subset  \{1,2,3,4\}$.

Conversely, for each $l \in{\IR}_{>0} \times {\IR}^5$, if $\phi_{ij}(l) \in (-1,1) $, then by the definition of $\phi_{1i}(l)$, and using (\ref{3-1edge-1})-(\ref{3-1angle-1}), we may show that ${\theta}_{ij}^1,\; {\theta}_{ik}^1,\;{\theta}_{jk}^1$ satisfy the triangular inequalities. Then there exists an unique vertex hyperbolic triangle $\triangle_1 $ having them as edge lengths. And there exists an unique vertex hyperbolic triangle $\triangle_2 $ having them as edge lengths in the same way.

Taking ${\alpha}_{1i} = {\cos}^{-1} \phi_{1i}(l) \in (0,\pi)$,  ${\alpha}_{2j} = {\cos}^{-1} \phi_{2j}(l) \in (0,\pi)$, we see that $\alpha_{12},\, \alpha_{13},\, \alpha_{14}$ are the inner angles of the vertex triangle $\triangle_1 $, and $\alpha_{21},\, \alpha_{23},\, \alpha_{24}$ are the inner angles of the vertex triangle $\triangle_2 $, hence they satisfy $\alpha_{12} = \alpha_{21}$, $ \alpha_{12}+ \alpha_{13}+ \alpha_{14} < \pi $ and $ \alpha_{21}+ \alpha_{23}+ \alpha_{24} < \pi $. Since  $\alpha_{34} = \pi- \alpha_{13}+\alpha_{23}$, then all the six dihedral angles $ ( \alpha_{12}, \ldots, \alpha_{34}) $ are determined, and by directly calculations, their cosine values are also equal to the values given by Lemma \ref{prop 2.1}, hence we show that $l$ is the edge length vector of a $2$-$2$ type tetrahedron.
\end{proof}

\subsection{The geometric structure of $\CL$ and $\CR$}
Recall the definition of $\phi_{pq}$ for $2$-$2$ type and $3$-$1$ type hyperbolic tetrahedra, and $\CR$ respectively equals to ${\IR}_{>0} \times {\IR}^5$ or  ${\IR}^2_{>0} \times {\IR} \times {\IR}_{>0} \times {\IR}^2$ for $2$-$2$ type or $3$-$1$ type hyperbolic tetrahedra case respectively. For $\{i,j\} \subset  \{1,2,3,4\}$, let
$$\Omega_{ij}^{\pm} = \{l=(l_{12},\ldots,l_{34})\in \CR \mid  \pm \phi_{ij}(l) \geq 1 \},$$
and let
$$X_{ij}^{\pm} = \{(l_{12},\ldots,l_{34})\in \CR \mid   \phi_{ij}(l) = \pm 1 \},$$
then $\Omega_{ij}^{\pm}=\Omega_{ji}^{\pm}$, $X_{ij}^{\pm}=X_{ji}^{\pm}$ and $X_{ij}^{\pm}=\partial\Omega_{ij}^{\pm}$. Moreover, by Proposition~\ref{L-space}, for both of these two cases, we have
\begin{equation}
\CR \setminus{\CL}=\underset{i\neq j}{\cup} ({\Omega}_{ij}^{+} \cup {\Omega}_{ij}^{-}).
\end{equation}

\begin{lemma}\label{Set-relation}
In both the two cases, let $\phi$ be defined in the Definition \ref{def-phi-22} and Definition \ref{def-phi-31}, and correspondingly define $\Omega_{ij}^{\pm}$, $X_{ij}^{\pm}$ as above. For $\{ i, j, k, h \}=\{1,2,3,4\}$, we have
\begin{enumerate}
    \item[(1)] ${\Omega}_{ij}^{-}={\Omega}_{kh}^{-}$ and ${\Omega}_{ij}^{+}={\Omega}_{kh}^{+}$.\\[-4pt]
    \item[(2)] $ {\Omega}_{ij}^{-}\cap{\Omega}_{ik}^{-}=\emptyset$, ${\Omega}_{ij}^{+}={\Omega}_{ik}^{-}\sqcup{\Omega}_{ih}^{-}$ and
    ${\Omega}_{ij}^{-}={\Omega}_{ik}^{+}\cap{\Omega}_{ih}^{+}$.\\[-4pt]
    \item[(3)] ${X}_{ij}^{-} = {X}_{kh}^{-}$, ${X}_{ij}^{+}={X}_{kh}^{+}$, ${X}_{ij}^{+}={X}_{ik}^{-}\sqcup{X}_{ih}^{-}$ and ${X}_{ij}^{-}={X}_{ik}^{+}\cap{X}_{ih}^{+}$.\\[-4pt]
    \item[(4)] For each pair $\{i,j\}\subset\{1,2,3,4\}$, we have
    $$\CR \setminus{\CL}={\Omega}_{12}^{-}\sqcup{\Omega}_{13}^{-}\sqcup{\Omega}_{14}^{-}={\Omega}_{ij}^{-} \sqcup {\Omega}_{ij}^{+}.$$
\end{enumerate}
\end{lemma}
\begin{proof}
The proof strategies for $2$-$2$ type case and $3$-$1$ type case are similar. As an example, here we only sketch the proof of conclusion (1) in $2$-$2$ type case, and the proofs of remaining conclusions could be obtain using the same argument of Lemma 2.10 in \cite{Feng-Ge-Liu}. Notice that for each $l\in \CR$, the following $\theta^i_{jk}(l) (\{i,j,k\}\subset\{1,2,3,4\})$ all make sense by Theorem \ref{thm-non-degenerate}.

~

For $2$-$2$ type tetrahedra case, let
\[x = e^{l_{13} + l_{14} - l_{34}}, y = e^{l_{14} + l_{23} -l_{34}}, z = e^{l_{13} + l_{24} - l_{34}}, d = \cosh {l_{12}},\]
then
\[\phi_{12} = \frac{(x^2 + xyd + xzd + yzd^2) - (2x + yz)(d^2 - 1)}{\sqrt {x^2y^2 + y^2z^z + x^4 + x^2z^2 + 2xy^zd + 2xyz^2d + 2x^3yd +2x^3zd +4x^2yzd^2}},\]
\[\phi_{34} = \frac{xy + xz -2xd}{2\sqrt {(yz + x)(x^2 + x)}}.\]
By direct calculations, it's easy to check that  $\pm\phi_{12}(l)\geq 1$ if and only if $\pm\phi_{34}(l)\geq 1$. Hence ${\Omega}_{12}^{-}={\Omega}_{34}^{-}$ and ${\Omega}_{12}^{+}={\Omega}_{34}^{+}$. Similarly, ${\Omega}_{13}^{-}={\Omega}_{24}^{-}$, ${\Omega}_{13}^{+}={\Omega}_{24}^{+}$, ${\Omega}_{14}^{-}={\Omega}_{23}^{-}$ and ${\Omega}_{14}^{+}={\Omega}_{23}^{+}$. Hence conclusion (1) is obtained.
\end{proof}

Following the notation above,
let $ {\Omega}_{1} = {\Omega}_{12}^{-}$, \; ${\Omega}_{2} = {\Omega}_{13}^{-}$ and $ {\Omega}_{3} = {\Omega}_{14}^{-}$.
Similarly, let $ {X}_{1} = {X}_{12}^{-}$, \; ${X}_{2} = {X}_{13}^{-}$ and $ {X}_{3} = {X}_{14}^{-}$.
By Lemma \ref{Set-relation}, we know that $X_i = \partial {\Omega}_{i}$ and
 \[
   \partial \CL \subset  \partial {\Omega}_{1} \cup \partial {\Omega}_{2} \cup \partial {\Omega}_{3} = {\sqcup}_{i=1}^3 X_{i} .
 \]

\begin{prop}\label{L-frontier-22}
In both the two cases, let $\partial \CL$ be the frontier of $\CL$ in $\CR$, then
  $X_{1}$, $X_{2}$, $X_{3}$ are
  real analytic codimension-1 submanifold of $\CR$ and
\[ \partial \CL = {\sqcup}_{i=1}^3 X_{i}.\]
The complement $\CR\setminus \CL$ is a disjoint union of three analytic manifolds ${\Omega}_{i} $ with boundary so that ${\Omega}_{i} \cap \partial \CL = X_i, i = 1, 2, 3.$
\end{prop}

\begin{proof}
We could get this proposition by similar argument in Proposition 2.11 in \cite{Feng-Ge-Liu} and Proposition 4.5 in \cite{Luo2018}.

For each $(l_{12},\cdots,l_{34}) \in \CR$, let $Q_{ij}$ be the denominator of ${\phi}_{ij}(l)$.

First consider the $2$-$2$ case. In this case, recall $l_{ij} = l_{ji}$, we have
  \[\frac{\partial {\phi}_{12}(l)}{\partial l_{34}} = -\frac{2e^{l_{34}}(\sinh l_{12})^2}{Q_{12}},\]
  \[\frac{\partial {\phi}_{34}(l)}{\partial l_{12}} = -\frac{2e^{l_{34}}\sinh l_{12}}{Q_{34}},\]
and other partial derivatives:
  \[\frac{\partial {\phi}_{ik}(l)}{\partial l_{jh}} = -\frac{e^{l_{jh}+l_{ik}-l_{ih}}}{ Q_{ij}}.\]

It's clear that these are not zero for any $(l_{12},\ldots,l_{34}) \in \CR $,
  so $ \nabla \phi_{ij} (l) \neq 0$.

  By the Implicit Function Theorem, since $-1$ is a regular value of $\phi_{ij} (l)$, $X_i = {\phi_{ij}^{-1}(-1)} $ is a smooth
  co-dimension-1 submanifold of $\CR $.
  Further, since each $\phi_{ij}$ is real analytic in $\CR $,
  $X_{1}$, $X_{2}$ and $X_{3}$ are real analytic codimension-1 submanifold of $\CR$.

  By Lemma \ref{Set-relation},
  it's clear that $\partial \CL \subset {\sqcup}_{i=1}^3 X_{i}.$ Now we claim that $ X_i \subset  \partial \CL $. Without loss of generality, here we just show $X_1 \subset  \partial \CL$.
  Indeed, for each $l= (l_{12},\ldots,l_{34}) \in X_1 $,
  let $\epsilon_n \to 0^{+}$, define $l^{(n)} = (l_{12}- \epsilon_n, l_{13}, l_{14}, l_{23},l_{24}, l_{34}) \to l $. By the definition of $X_1 $ and calculating the partial derivative $\frac{\partial \theta_{13}^2}{\partial l_{14}}$, $\frac{\partial \theta_{14}^2}{\partial l_{14}}$, then we know that, for $n$ large enough,
  ${\theta_{13}^2}^{(n)} ,\, {\theta_{14}^2}^{(n)}, \,{\theta_{34}^2}^{(n)} $ satisfy the
  triangular inequalities, and $\phi_{12}(l^{(n)}),\phi_{23}(l^{(n)}),\phi_{24}(l^{(n)}) \in (-1,1)$. By similar calculation, each $\phi_{ij}(l^{(n)}) \in (-1, 1)$, thus $l^{(n)} \in \CL$ by Proposition \ref{L-space}. Thus, $l \in \partial \CL$, this also means $X_1 \subset  \partial \CL$. As this point, we have completed the proof of the claim.
  Therefore, $ \partial \CL = {\sqcup}_{i=1}^3 X_{i} $.

Next consider the $3$-$1$ case. In this case, let $\{i,j,k\}=\{1,2,3\}$, then we have

  \[\frac{\partial {\phi}_{i4}(l)}{\partial l_{jk}} = -\frac{e^{2l_{i4}}\sinh l_{jk}}{Q_{i4}} \]
   and others:
  \[\frac{\partial {\phi}_{ij}(l)}{\partial l_{k4}} = -\frac{e^{l_{k4}}\sinh^2 l_{ij}}{Q_{ij}}.\]
  It's clear that these are not zero for any $(l_{12},\ldots,l_{34}) \in \CR$,
  so $ \nabla \phi_{ij} (l) \neq 0$. Therefore, $-1$ is a regular value of $\phi_{ij} (l)$.
  Then we may finish the remainder of this proof by using the same argument as in the $2$-$2$ case.
\end{proof}

\subsection{Extending $\alpha$ from $\CL$ to $\overline{\CR}$ continuously}
\label{section-extend-a-to-bar-R}
In this subsection, we extend the definition of those dihedral angles $\{\alpha_{pq}\}$ from $\CL$ to $\overline{\CR}$ continuously following the spirit of \cite{Bobenko2015,Luo2011,Luo2018}.

In the $2$-$2$ type case, the extension is direct. Let $\overline{\CL}$ be the closure of $\CL$ in $\overline{\CR}={\IR}_{\ge 0}\times {\IR}^5$, for each pair $\{i,j\}\subset\{1,2,3,4\}$, define $a_{ij}\vert_{ \overline{\Omega_{ij}^{+}}}=0$ and
$a_{ij}\vert_{ \overline{\Omega_{ij}^{-}}}=\pi$. As a consequence of Proposition \ref{L-space}, Lemma \ref{Set-relation} and Proposition \ref{L-frontier-22},
we have

\begin{cor}\label{angle-extended-22}
In a $2$-$2$ type hyperbolic tetrahedron, for each pair $\{i,j\}\subset\{1,2,3,4\}$, the dihedral angle $\alpha_{ij}: \CL \to \IR$ can be extended continuously to $\overline{\CR}$, so that its extension,
still denoted by $\alpha_{ij}: \overline{\CR}\to \IR $, is a constant on each component of $\overline{\CR}\setminus \overline{\CL}$.
\end{cor}

For the $3$-$1$ type case, the situation is slightly more complicated. Recall in this case, $\overline{\CL}$ is the closure of $\CL$ in $\overline{\CR}={\IR}^2_{\ge 0} \times {\IR} \times {\IR}_{\ge 0} \times {\IR}^2$. Let $\overline X_i$ be the closure of $X_i$ in $\overline{\CR}$.

\begin{definition}
\label{subset-31}
In a $3$-$1$ type tetrahedron, let $S$ be any subset of the edges $\{e_{ij}\}_{1\leq i<j\leq 4}$ which includes $\{e_{14}, e_{24}, e_{34}\}$. Define $\CD_S$ to be the space composed of all generalized metrics $l \in \overline{\CR}$, so that $l_{ij}>0$ for $e_{ij}\in S$ and $l_{ij}=0$ for $e_{ij}\notin S$, where $\{i,j\}\subset\{1,2,3\}$.
\end{definition}

By definition, we see $\overline{\CR}=\sqcup_S\CD_S$, where $S$ runs over every subset of the edges $\{e_{ij}\}_{1\leq i<j\leq 4}$ which contains $\{e_{14}, e_{24}, e_{34}\}$.

\begin{prop}\label{subset-mfd-31}
If $\CD_S \cap \overline X_i \neq \emptyset$, then $X_i^S \doteq \CD_S \cap \overline X_i$ is a real analytic codimension-1 submanifold of $\CD_S$.
\end{prop}
\begin{proof}
For $\{i,j\}\subset\{1, 2, 3\}$, let $\overline {X^-_{ij}}$ be the closure of $X^-_{ij}$ in $\overline{\CR}$.
If $e_{ij} \notin S$, then $l_{ij}=0$. By Lemma \ref{L-zero-31}, $\phi_{ij}(l) = 1$. On the other hand, if $l \in \overline {X^-_{ij}}$, then $\phi_{ij}(l) = -1$. Thus, if $e_{ij} \notin S$, then $\CD_S \cap \overline {X^-_{ij}} = \emptyset$.
So, if $l \in \CD_S \cap \overline {X^-_{ij}}$, then $e_{ij} \in S$. Thus $l_{ij} > 0$ and $\sinh l_{ij} \neq 0$. Further, we have
\[\frac{\partial \phi_{ij}}{\partial l_{k4}} = -\frac{e^{l_{k4}}\sinh^2 l_{ij}}{Q_{ij}} \neq 0,\]
which implies that the projection of $\nabla \phi_{ij}$ to tangent space of $\CD_S$ at $l\in \CD_S \cap \overline {X^-_{ij}}$ is non-vanishing, this means that $\CD_S$ and $\overline {X^-_{ij}}$ transversely intersect. By the Implicit Function Theorem, the intersection is a smooth codimension-1 submanifold of $\CD_S$. Since each $\phi_{ij}$ is real analytic, the submanifold is real analytic in $\CD_S$.
\end{proof}
Set $\CL_S = \CD_S \cap \overline {\CL}$. For each pair $\{i,j\}\subset\{1,2,3,4\}$, define $a_{ij}\vert_{ \overline{\Omega_{ij}^{+}}}=0$ and
$a_{ij}\vert_{ \overline{\Omega_{ij}^{-}}}=\pi$. By Proposition \ref{L-space}, Lemma \ref{Set-relation}, Proposition \ref{L-frontier-22} and Proposition \ref{subset-mfd-31}, we have

\begin{cor}\label{angle-extended-31}
In a $3$-$1$ type hyperbolic tetrahedron, for each pair $\{i,j\}\subset\{1,2,3,4\}$, the dihedral angle $\alpha_{ij}: \CL \to \IR$ can be extended continuously to $\overline{\CR}$, so that its extension, still denoted by $\alpha_{ij}$, is a constant on each component of $\CD_S \setminus \CL_S$ for each set $S$ in Definition \ref{subset-31}.
\end{cor}

\section{The volume and co-volume function of hyperbolic tetrahedra}
\label{section-volume}
In this section, we will introduce the Schl\"{a}fli formula, as well as the volume and co-volume functions and their extensions to general hyperbolic tetrahedra.

\subsection{The volume and co-volume function}
For an ideal tetrahedron, its opposite sides have the same dihedral angle. Assuming that the dihedral angles on its three sets of opposite sides are $\theta_1$, $\theta_2$ and $\theta_3$, respectively, then its volume is $\Lambda(\theta_1)+\Lambda(\theta_2)+\Lambda(\theta_3)$, where $\Lambda(\theta)$ is Milnor's Lobachevsky function defined by
\[\Lambda(\theta) = -\int_{0}^{\theta} \ln \vert 2\sin t \vert \di t.\]
which is $\pi$-periodic, odd and smooth except at $\theta\in \pi\mathbb{Z}$, with tangents vertical (\cite{MR1277810}). Note $\theta_1+\theta_2+\theta_3=\pi$ in an ideal tetrahedron, and Rivin \cite{Rivin1994} proved that the volume is strictly concave down on $\{(\theta_1, \theta_2, \theta_3)\in(0, \pi)^3:\theta_1+\theta_2+\theta_3=\pi\}$.

\begin{figure}[htbp]
	\centering
	{\includegraphics[width=.50\columnwidth]{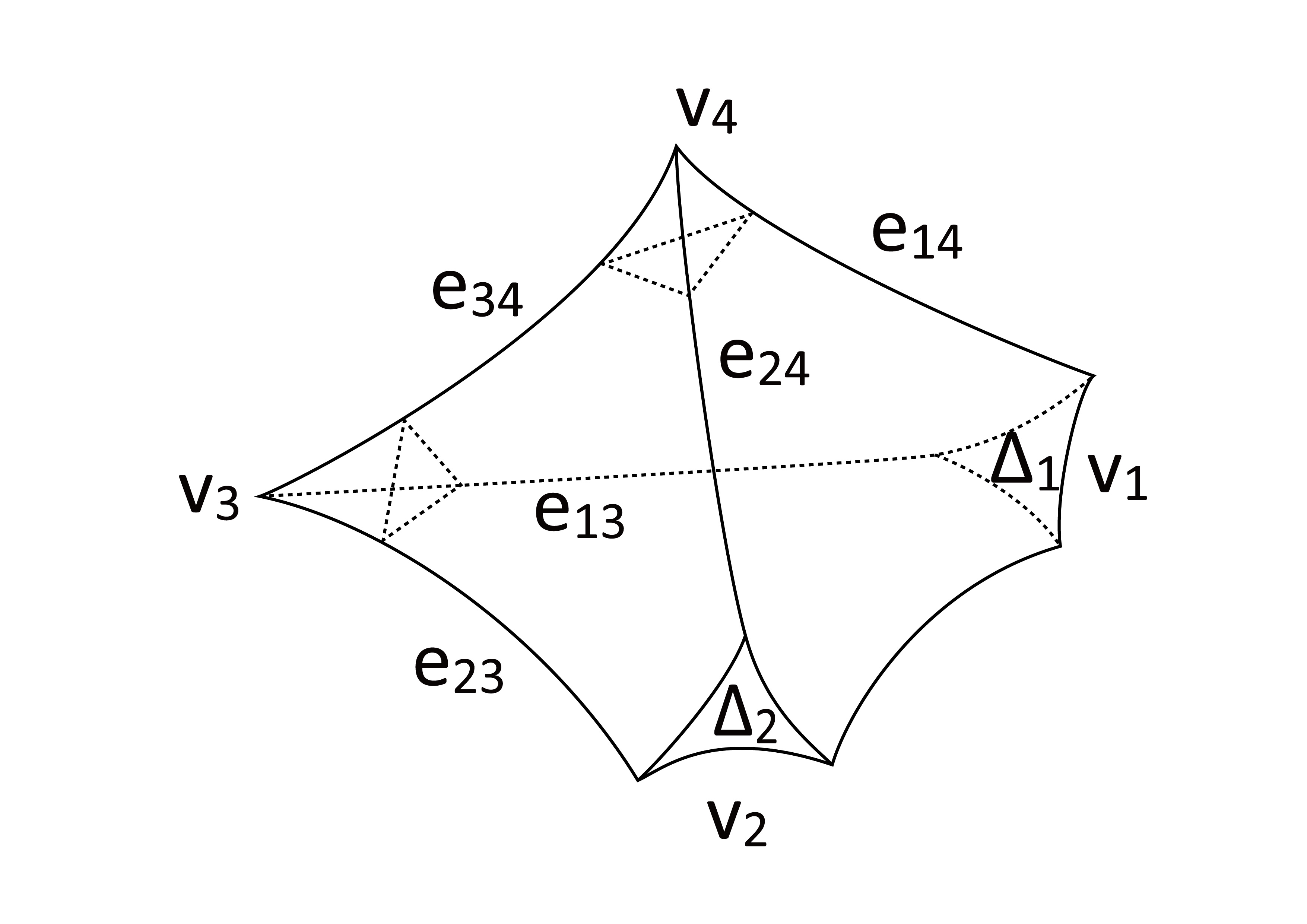}}
	\caption{a $2$-$2$ type hyperbolic tetrahedron}
    \label{figure4}
\end{figure}

Now we consider a $2$-$2$ type hyperbolic tetrahedron $T_{2}$. Recall our notation conventions, the two hyper-ideal vertices are denoted respectively by $v_1$ and $v_2$, the two ideal vertices are denoted respectively by $v_3$ and $v_4$, the dihedral angle at an edge $e_{ij}$ is denoted by $\alpha_{ij}$ for each pair $\{i,j\}\subset\{1,2,3,4\} $, see Figure \ref{figure4}. Then the volume of $T_{2}$ is determined by its six dihedral angles, that is,
\begin{eqnarray*}
   2 vol(T_{2}) &=& \Lambda(\alpha_{13}) + \Lambda(\alpha_{23}) + \Lambda(\alpha_{14}) + \Lambda(\alpha_{24})\\
              &+&\Lambda \bigl(  \frac{\pi-\alpha_{12}-\alpha_{13}+\alpha_{14}}{2}  \bigr)
              +\Lambda \bigl(  \frac{\pi-\alpha_{13}-\alpha_{14}+\alpha_{12}}{2}  \bigr)
              +\Lambda \bigl(  \frac{\pi-\alpha_{12}-\alpha_{13}-\alpha_{14}}{2}  \bigr)\\
              &+&\Lambda \bigl(  \frac{\pi-\alpha_{23}-\alpha_{24}+\alpha_{12}}{2}  \bigr)
              +\Lambda \bigl(  \frac{\pi-\alpha_{12}-\alpha_{24}+\alpha_{23}}{2}  \bigr)
              +\Lambda \bigl(  \frac{\pi-\alpha_{24}-\alpha_{23}-\alpha_{12}}{2}  \bigr)\\
              &-&\Lambda \bigl(  \frac{\pi+\alpha_{13}-\alpha_{14}-\alpha_{12}}{2}  \bigr)
              -\Lambda \bigl(  \frac{\pi+\alpha_{24}-\alpha_{23}-\alpha_{12}}{2}  \bigr).
\end{eqnarray*}

Next we consider a $3$-$1$ type hyperbolic tetrahedron $T_{3}$. In this case, $v_1$, $v_2$ and $v_3$ are the three hyper-ideal vertices, and $v_4$ is the ideal vertex. Then the volume of $T_{3}$ is also determined by its six dihedral angles, that is,
\begin{eqnarray*}
   2 vol(T_{3}) &=& \Lambda(\alpha_{14}) + \Lambda(\alpha_{24}) + \Lambda(\alpha_{34}) \\
              &+&\Lambda \bigl(  \frac{\pi-\alpha_{12}-\alpha_{23}+\alpha_{24}}{2}  \bigr)
              +\Lambda \bigl(  \frac{\pi-\alpha_{12}-\alpha_{24}+\alpha_{23}}{2}  \bigr)\\
              &+&\Lambda \bigl(  \frac{\pi-\alpha_{12}-\alpha_{23}-\alpha_{24}}{2}  \bigr)
              -\Lambda \bigl(  \frac{\pi+\alpha_{12}-\alpha_{23}-\alpha_{24}}{2}  \bigr)\\
              &+&\Lambda \bigl(  \frac{\pi-\alpha_{34}-\alpha_{23}+\alpha_{13}}{2}  \bigr)
              +\Lambda \bigl(  \frac{\pi-\alpha_{13}-\alpha_{34}+\alpha_{23}}{2}  \bigr)\\
              &+&\Lambda \bigl(  \frac{\pi-\alpha_{13}-\alpha_{23}-\alpha_{34}}{2}  \bigr)
              -\Lambda \bigl(  \frac{\pi+\alpha_{34}-\alpha_{23}-\alpha_{13}}{2}  \bigr)\\
              &+&\Lambda \bigl(  \frac{\pi-\alpha_{12}-\alpha_{14}+\alpha_{13}}{2}  \bigr)
              +\Lambda \bigl(  \frac{\pi-\alpha_{12}-\alpha_{13}+\alpha_{14}}{2}  \bigr)\\
              &+&\Lambda \bigl(  \frac{\pi-\alpha_{12}-\alpha_{13}-\alpha_{14}}{2}  \bigr)
              -\Lambda \bigl(  \frac{\pi+\alpha_{12}-\alpha_{13}-\alpha_{14}}{2}  \bigr).
\end{eqnarray*}

Consider the volume function $vol(\cdot)$, defined on $\CB_{2|2}$ or $\CB_{3|1}$ respectively, the Schl\"{a}fli formula \cite{Bonahon1998,Schlenker2002} have a consistent form for both $2$-$2$ type and $3$-$1$ type tetrahedra. It says
\begin{equation}\label{equ-shlafi}
\di vol=-\frac{1}{2}\sum_{1\leq i<j\leq 4}l_{ij}d\alpha_{ij}.
\end{equation}

Let $\overline {\CB}_{2|2}$ (or $\overline {\CB}_{3|1}$, respectively) be the closure of $\CB_{2|2}$ (or $\CB_{3|1}$, respectively) in $\IR^6$, i.e.
\[
  \overline {\CB}_{2|2} =\bigl\{  (a_{12},\ldots,a_{34}) \in {\IR}^6_{\geq 0} \mid \sum_{j\neq k} {\alpha_{kj}}\le \pi,  k \in \{1,2\} , \sum_{j\neq i} a_{ij} = \pi \; \text{for}\; i\in \{3,4\}\bigr\},
\]
and
\[
  \overline {\CB}_{3|1} =\bigl\{  (a_{12},\ldots,a_{34}) \in {\IR}^6_{\geq 0} \mid \sum_{j\neq k} {\alpha_{kj}}\le \pi,  k \in \{1,2,3\} , \sum_{j\neq 4} a_{j4} = \pi  \bigr\}.
\]

By Rivin \cite{Rivin2008}, the volume function can be extended continuously to the compact closure $\overline{\CB}$ of $\CB$. Leibon~\cite{Leibon2002} and Schlenker~\cite{Schlenker2002} proved the following lemma.

\begin{lemma}[\cite{Leibon2002} ~\cite{Schlenker2002} ]
Let $vol$ be the volume function of a $2$-$2$ type (or $3$-$1$ type, respectively) hyperbolic tetrahedron. Then the Hessian matrix of $vol$ with respect to dihedral angles is negative definite, therefore the volume function $vol$ is concave on $\overline {\CB}_{2|2}$ (or $\overline {\CB}_{3|1}$, respectively).
\end{lemma}

Recall $\CL$ is the space of vectors $(l_{12}, \ldots, l_{34})$ such that there exists a  $2$-$2$ type (or $3$-$1$ type, respectively) hyperbolic tetrahedron having $l_{ij}$ as the decorated length of edge $e_{ij}$. By Lemma \ref{angleedge31} and Remark \ref{unchange31}, if we set $ v_1 = (0,1,0,1,0,1)$ and $v_2 = (0,0,1,0,1,1)$, then $\CL$, $\alpha$, $vol$ are all invariant along the direction $v_1$, $v_2$ (or $v_2$, respectively) in the case of $2$-$2$ type (or $3$-$1$ type, respectively) tetrahedra. Denote $L(v_1,v_2)=\text{span}\{v_1, v_2 \}$ and $L(v_2)=\text{span}\{ v_2 \}$. From the definition of $\CB_{2|2}$ and $\CB_{3|1}$, they are convex and relative open subsets of $\IR^6$ and hence simply connected. By Proposition \ref{prop 2.1}, $\CB_{2|2}$ is diffeomorphic to $\CL / L(v_1,v_2)$, and $\CB_{3|1}$ is diffeomorphic to $\CL / L(v_2)$.
Thus both $\CB_{2|2}$ and $\CL / L(v_1,v_2)$ can be considered as simply connected open sets in $\mathbb{R}^4$, and both $\CB_{3|1}$ and $\CL / L(v_2)$ can be considered as simply connected open sets in $\mathbb{R}^5$. Consequently, by Remark \ref{unchange31}, we see $\CL$ is a simply connected open subset of $\IR^6$. Moreover, the volume function $vol(\cdot)$ may be considered as a function defined on $\CL$.

Now, for the sake of convenience and without causing confusion, we will abbreviate the set of dihedral angles $\CB_{2|2}$ and $\CB_{3|1}$ in different cases as $\CB$.

\begin{definition}[co-volume]
The co-volume function $ cov(l):\CL \rightarrow \IR$ is exactly the Fenchel-Legendre dual transform of $vol(\alpha):{\CB} \rightarrow \IR$, and takes the following form
\begin{equation}
  cov(l) = 2vol(l) + \sum_{i<j} \alpha_{ij}\cdot l_{ij} .
\end{equation}
\end{definition}

\begin{thm}\label{cov-po}

The co-volume function $cov:\CL \rightarrow \IR$ is locally convex. To be precise,
\begin{enumerate}
  \item for the case $\CB=\CB_{2|2}$, the co-volume function has a semi-positive definite Hessian matrix with rank 4 at each $l\in \CL$. Moreover, $cov$ is locally strictly convex on the open set $\CL / L(v_1,v_2)\subset \mathbb{R}^4$.
  \item for the case $\CB=\CB_{3|1}$, the co-volume function has a semi-positive definite Hessian matrix with rank 5 at each $l\in \CL$. Moreover, $cov$ is locally strictly convex on the open set $\CL / L(v_2)\subset \mathbb{R}^5$.
\end{enumerate}
\end{thm}
\begin{proof}
   We consider $cov$ as a function of the edge length variables $(l_1,\ldots,l_6)$, and
   the dihedral angle at edge $l_{j}$ is $ \alpha_{j}$, by the Schl\"{a}fli formula, then
$$\di cov = 2 \di vol + \sum_{i=1}^{6}(\alpha_{i} \di l_{i} + l_{i}\di \alpha_{i}) = \sum_{i=1}^{6} \alpha_{i} \di l_{i}.$$
Thus $\partial cov/\partial l_j=\alpha_j$ for each $1\leq j\leq 6$. As a consequence, the Hessian of $cov$ with respect to $(l_1,\ldots,l_6)$ is exactly ${[\partial  \alpha_j / \partial l_i]}_{6\times 6}$. The conclusion may be derived by the similar argument as in the proof of Theorem 2.6 in \cite{Feng-Ge-Liu}, and we leave out the details here.
\end{proof}

From the Schl\"{a}fli formula (\ref{equ-shlafi}) and the proof of Theorem \ref{cov-po}, for both $2$-$2$ type and $1$-$3$ type hyperbolic tetrahedra, the locally convex function co-volume function
$ cov : \CL \to \IR$ satisfies
\[
  \frac{\partial cov}{\partial l_{ij}} = {\alpha}_{ij}
\]
for each pair $\{i,j\}\subset\{1,2,3,4\}$, where ${\alpha}_{ij} : \CL \to \IR $ is the dihedral angle function. In particular, the differential $1$-form $w=\sum_{1\leq i<j\leq4} {\alpha}_{ij} d l_{ij}$ is smooth and closed in the simply connected open set $\CL\subset \IR^6$, and we can recover $cov$ by the integration $cov(l) = \int^{l} w$.

\subsection{The extended co-volume function}
Now its time to extend the definition of those dihedral angles $\{\alpha_{ij}(l)\}$, by extending their domain of definition from $l\in\CR$ to ${\IR}^6$.
\begin{definition}[edge length vector $l^+$]
\label{def-l+}
For each $l=(l_{12},\cdots,,l_{34}) \in {\IR}^6$, if $e_{ij}$ is an edge with decorated edge length $l_{ij} < 0$, then endow it a new length $l^{+}_{ij} = 0$, otherwise, $l_{ij}$ doesn't change. In this way, we get a new edge length vector $l^+$.
\end{definition}

For example, in a $3$-$1$ type hyperbolic tetrahedron, each generalized metric $l\in\IR^6$ corresponds to a new metric
$$l^+=(l_{12}^+,l_{13}^+,l_{14},l_{23}^+,l_{24},l_{34})\in\overline{\CR}={\IR}^2_{\ge0}\times{\IR}\times{\IR}_{\ge0}\times{\IR}^2,$$
where $l_{ij}^+ = \max {\{0, l_{ij}\}}$. Hence the dihedral angle $\alpha_{ij}: \CR \to \IR $ is extended continuously to
$$\alpha_{ij}: {\IR}^6 \to \IR$$
by setting $\alpha_{ij}(l)=\alpha_{ij}(l^+)$. The extended functions are still denoted by $\alpha_{ij}$.

\begin{definition}[generalized hyperbolic tetrahedron]
A generalized $2$-$2$ type (or $3$-$1$ type, respectively) hyperbolic tetrahedron $\sigma$ is a purely topological and combinatorial tetrahedron that has the same combinatorial structure as a normal $2$-$2$ type (or $3$-$1$ type, respectively) hyperbolic tetrahedron, with each edge $e_{ij}$ assigned a real number $l_{ij}$ as the edge length.
\end{definition}


Knowing which type of tetrahedron it is, each $l\in \IR^6$ determines a generalized hyperbolic tetrahedron. Based on Corollary \ref{angle-extended-22}, Corollary \ref{angle-extended-31} and Definition \ref{def-l+} above, the dihedral angle on each edge $e_{ij}$ of the generalized $2$-$2$ type or $3$-$1$ type hyperbolic ideal tetrahedron could be obtained by $\alpha_{ij}(l)$. For general $1$-$3$ type, $0$-$4$ type and $4$-$0$ type hyperbolic tetrahedra, corresponding definitions are given in \cite{Feng-Ge-Liu} and \cite{Luo2018}. Now we are going to define the co-volume function of a generalized $2$-$2$ or $3$-$1$ type hyperbolic tetrahedron in the following of this section. Consider a new continuous differential $1$-form $\mu$ on ${\IR}^6$ given by
\[
  \mu(l) = \sum_{i\neq j} \alpha_{ij}(l) \di l_{ij} .
\]

\begin{lemma}[\cite{Luo2011}]\label{Luo-lemma}
Suppose $U \subset  \IR^N$ is an open set and $\mu(x)=\sum_i{\alpha}_{i}(x) \di x_{i}$ is a continuous differential 1-form on $U$ (by definition, $\mu$ is called continuous if all its coefficients $\alpha_i(x)$ are continuous. Moreover, $\mu$ is called closed if $\int_{\partial \Delta} \mu = 0$ for any Euclidean triangle $\Delta\subset U$).
  \begin{enumerate}
    \item If $ A \subset  U $ is an open subset bounded by a smooth codimension-1 submanifold of $U$, both $\mu\vert_A$ and $\mu\vert_{U\setminus \bar{A}} $ are closed,
    then $\mu$ is closed in $ U $.
    \item If $U$ is simply connected, then $F(x) =\int^x \mu $ is a $C^1$-smooth function such that \[ \frac{\partial F}{\partial x_i} = \alpha_i . \]
    \item If $U$ is convex and $A \subset  U$ is an open subset of $U$ bounded by a codimension-1 real analytic submanifold of $U$ so that $ F \vert_A$ and $F \vert_{U\setminus \bar{A}}$ are locally convex, then $F$ is convex in $U$.
  \end{enumerate}
\end{lemma}
\begin{prop}\label{closed-1-form}
The continuous differential 1-form $\mu(l)=\sum_{i\neq j} {\alpha}_{ij}(l) \di l_{ij}$ is closed in $\IR^6$, that is, for any Euclidean triangle $\Delta$ in $\IR^6$, $\int_{\partial \Delta} \mu = 0$.
\end{prop}

\begin{proof}
For $2$-$2$ type and $3$-$1$ type hyperbolic tetrahedra cases, we give the proof respectively.

(1) The $2$-$2$ type tetrahedra case. Take $U=\CR={\IR}_{> 0}\times {\IR}^5$ and $A=\CL$. Obviously, ${\alpha}_{ij}(l) $ is continuous on $U$ by Corollary \ref{angle-extended-22}. On $A$, $\mu\vert_A=\omega$ is smooth and closed. On each connect component of component of $U\setminus A$, $\alpha_{ij}$ is constant by Corollary \ref{angle-extended-22}, and hence $\mu\vert_{U\setminus \bar{A}}$ is also closed. Then using part (1) of Lemma \ref{Luo-lemma}, $\mu$ is closed in $\CR$. For any triangle $\Delta\subset\bar{\CR}$, it can always be approximated by triangles in $\CR$. Further by the continuity of $\mu$ derived in $\int_{\partial \Delta} \mu = 0$, implying that $\mu$ is closed in $\bar{\CR}$.

Next we show that $\mu$ is closed in ${\IR}^6$. For any Euclidean triangle $\Delta\subset{\IR}^6 \setminus {\bar{\CR}}$, let $\Delta_{pr}$ be the projection of $\Delta$ to $\bar{\CR}$, by replacing $l$ to $l^+$. By the definition of $l^{+}$ and Corollary \ref{angle-extended-22}, we have $\alpha_{12}(l)=\alpha_{12}(l^+)=0$ if $l_{12} \le 0$, which implies that $\int_{\Delta}\mu=\int_{\Delta_{pr}}\mu =0$. Consequently, $\mu$ is closed in ${\IR}^6 \setminus {\bar{\CR}}$. Taking $U=\IR^6$ and $A=\CR$, and applying part (1) of Lemma \ref{Luo-lemma} again, we conclude that $\mu$ is closed in ${\IR}^6$.

(2) The $3$-$1$ type tetrahedra case. Similar to the first paragraph of (1), replacing Corollary \ref{angle-extended-22} with Corollary \ref{angle-extended-31}, we show that $\mu$ is closed in both
We further prove that $\mu$ is closed in ${\IR}^6$ as follows. According to Definition \ref{def-phi-31}, Lemma \ref{angleedge31}, Lemma \ref{L-zero-31} and Definition \ref{subset-31}, we can write the concrete form of $\alpha_{ij}(l)$ and hence $\mu(l)$ for $l\in\CL_S$, where $S$ is any subset of the six edges (at lest includes $\{e_{14}, e_{24}, e_{34}\}$). Tedious but straightforward calculations show that $\mu|_{\CD_S}$ is closed in $\CL_S$ (one may give an alternate proof by approximating a triangle in $\CL_S$ by triangles in $\CL$). In addition, $\mu|_{\CD_S}$ is constant in each connected component of $\CD_S \setminus \CL_S$ by Corollary \ref{angle-extended-31}. Since Proposition \ref{subset-mfd-31} shows that the subset $\CL_S$ in $\CD_S$ is open and bounded by a smooth codimension-1 submanifold, we could use part (1) of Lemma \ref{Luo-lemma} for $U=\CD_S$ and $A=\CL_S$, to show that $\mu|_{\CD_S}$ is closed in ${\CD_S}$.

Next, we call the domain $Q = Q_1 \times Q_2 \times \IR \times Q_3 \times \IR^2$ a ``big quadrant" in ${\IR}^6$, where $Q_1, Q_2$ and $Q_3$ are respectively any open quadrant of $\IR$. For each Euclidean triangle $\Delta$ in a ``big quadrant" $ Q$ of $\IR^6$, and let $\Delta_S$ be the projection of $\Delta$ to ${\CD_S}$. By the definition of $\alpha_{ij}(l)$ and Lemma \ref{angle-extended-31}, $\alpha_{ij}(l) = 0$ if $l_{ij} \le 0$ for $\{i,j\}\subset\{2,3,4\}$, which implies $\int_{\Delta}\mu=\int_{\Delta_S}\mu =0$. As a consequence, $\mu$ is closed in $ Q$.
Since $\mu$ is closed in each ``big quadrant", then $\mu$ is closed in
$\CR=\IR^2_{>0}\times\IR\times\IR_{>0}\times\IR^2$ and $(\IR^2_{>0}\times\IR^4)\setminus\overline{\CR}$. Thus $\mu$ is closed in $\IR^2_{> 0} \times \IR^4$ by Lemma \ref{Luo-lemma} (1). Similarly, $\mu$ is closed in $\IR_{> 0} \times \IR_{< 0} \times \IR^4$. Further, by using Lemma \ref{Luo-lemma} (1) again, we obtain that $\mu$ is closed in $\IR_{> 0} \times \IR^5$. Likewise, $\mu$ is closed in $\IR_{< 0} \times \IR^5$ by repeating the above process. Finally, let's use Lemma \ref{Luo-lemma} (1) one more time, we conclude that $\mu$ is closed in the whole space ${\IR}^6$.
\end{proof}

\begin{cor}\label{general-cov}
The extended co-volume function $cov : \IR^6 \to \IR $ defined by the line integral
\[cov(l) = \int_{(0,\cdots,0)}^{l} \mu + cov(0,\cdots,0) \]
is a $C^1$-smooth convex function.
\end{cor}

\begin{proof}
In Proposition~\ref{closed-1-form}, we have proved that $\mu$ is continuous and closed in ${\IR}^6$.
Since ${\IR}^6$ is open and simply connected, thus $cov$ is a $C^1$-smooth function in ${\IR}^6$ by Lemma \ref{Luo-lemma} (2). On the other hand, by Lemma~\ref{Luo-lemma} (3), if let $A = \CL$, it is easy to show that $cov$ is convex in $\CR$.

For $2$-$2$ type case, $cov$ is convex in $\CR = {\IR}_{> 0}\times {\IR}^5$, and we know that $cov$ is convex in $\overline \CR$ because of  the continuity. On the other hand, if $l_{12} \le 0$, $a_{12}(l)=0$. This means that $cov(l) = cov(l^+)$. Thus, we know that $cov$ is convex in ${\IR}^6 \setminus {\overline \CR}$. Using Lemma \ref{Luo-lemma} (3), we conclude that $cov$ is convex in ${\IR}^6$.

For $3$-$1$ type case, by the continuity, we know that  $cov$ is convex  in $\overline {\CR}$. On the other hand, if $l_{ij} \le 0$ for $\{i,j\}\subset\{1, 2, 3\}$, $a_{ij}(l) = 0$, so $cov(l) = cov(l^+)$. Thus, we know that $cov$ is convex in each ``big quadrant" $Q$. Using the similar proof strategy in Proposition \ref{closed-1-form}, by repeating use of Lemma \ref{Luo-lemma} (3), we conclude that $cov$ is convex in ${\IR}^6$.
\end{proof}

\begin{rmk}\label{cov-mixed}
Luo-Yang showed that Corollary \ref{general-cov} is valid for ideal hyperbolic tetrahedra and hyper-ideal hyperbolic tetrahedra (see Lemma 2.2 and Corollary 4.12 in \cite{Luo2018}). For $1$-$3$ type hyperbolic tetrahedra, the same result is obtained by the authors and Liu (see Corollary 2.15 in \cite{Feng-Ge-Liu}). In summary, for any general hyperbolic tetrahedron, the extended co-volume function $cov : \IR^6 \to \IR $ could be defined, and is $C^1$-smooth convex.
\end{rmk}

\section{Global rigidity of hyperbolic polyhedral 3-manifolds}
\label{section-rigidity}
Starting from this section, we will shift our focus from studying individual tetrahedra to studying polyhedral 3-manifolds.
\subsection{Decorated polyhedral metric and combinatorial Ricci curvature}
\label{section-metric-curvature}
Let $(M,\CT)$ be a triangulated compact pseudo 3-manifold, that is, let $\{\sigma_1,\cdots, \sigma_t\}$ be a finite collection of combinatorial tetrahedra, so that $(M,\CT)$ is the quotient space of $\sigma_1 \sqcup\cdots \sqcup \sigma_t$, which is a simplicial complex of the disjoint union of combinatorial tetrahedra, via affine isomorphisms pairing faces of those tetrahedra. For the sake of convenience we will always assume the triangulation $\CT$ is \emph{closed}, that is, the gluing of tetrahedra around each internal edge is cyclical, or equivalently, each face in $\CT$ must be pasted onto another face in $\CT$. The relevant statements and conclusions can be easily extended to the situation of non-closed triangulations, with the main difference being a slight difference in the definition of curvature.

Denote $E=\{e_1, \dots, e_m\}$ by the set of edges in $\CT$, where $m$ is the number of edges. To simplify the notation, we often write $E=\{1,2,\cdots, m\}$, that is, each edge $e_i$ is replaced by the index $i$. Denote $V$ by the set of vertices in $\CT$. For any edge $e\in E$ and tetrahedron $\sigma\in\CT$ we say that $e$ is \emph{incident} to $\sigma,$ denoting by $e\prec\sigma$, if the former is contained in the latter. We also use $\CT$ to represent the set of tetrahedra, which is not confusing.

\begin{definition}[metric]
\label{dhpm}
A \emph{decorated hyperbolic polyhedral metric} on $(M,\CT)$ is obtained by replacing each combinatorial tetrahedron in $\CT$ by a partially truncated and decorated hyperbolic tetrahedron and replacing the affine gluing homeomorphisms by isometries preserving the decoration, i.e. gluing these hyperbolic tetrahedra along codimension-$1$ faces.
\end{definition}

A decorated hyperbolic polyhedral metric is often called a \emph{decorated metric}, or \emph{metric} in short. A metric is the same as assigning a number to each edge $e\in E$ so that each tetrahedron $\sigma$ becomes a partially truncated and decorated hyperbolic tetrahedron with assigned numbers as edge lengths. Denote
$$l=(l_1,l_2,\cdots l_m):=(l(e_1),l(e_2),\cdots, l(e_m))$$
by the edge length vector of a decorated metric, and denote
$$\CL(M, \mathcal{T})\subset \IR^E$$
by the space of all decorated metrics on $(M, \mathcal{T})$ parameterized by the edge length vector $l$.

\begin{definition}[curvature]
\label{def-curvature}
Each $l\in\CL(M, \mathcal{T})$ determines a \emph{combinatorial Ricci curvature} $K\in \IR^E$, which is defined as
$$K_e(l)=2\pi-\sum_{e\prec\sigma} \alpha_{\sigma, e}(l)$$
for each edge $e\in E$, where $\alpha_{\sigma, e}$ is the dihedral angle on $e$ in a hyperbolic tetrahedron $\sigma$.
\end{definition}

Given any $l \in \IR^6$, we can always get a generalized hyperbolic tetrahedron with corresponding truncations and decorations, any vector $l \in \IR^{E}$ assigns a number to each edge $e\in E$ so that each tetrahedron $\sigma$ becomes a generalized hyperbolic tetrahedron with assigned numbers as decorated edge lengths. Hence we may generalize the definition of a metric.

\begin{definition}[generalized metric]
\label{gdhpm}
Each $l\in\IR^E$ is called a \emph{generalized decorated hyperbolic polyhedral metric} (or \emph{generalized metric}, in short) on $(M, \mathcal{T})$.
\end{definition}

\begin{definition}[generalized curvature]
\label{defi:pre}
Each $l\in \IR^{E}$ determines a \emph{generalized combinatorial Ricci curvature} $\wt{K}\in \IR^E$, which is defined as
$$\wt{K}_e(l)=2\pi-\sum_{e\prec\sigma} \alpha_{\sigma, e}(l^+)$$
for each edge $e\in E$, where $\alpha_{\sigma, e}$ is the extended dihedral angle (see Section \ref{section-extend-a-to-bar-R}) on $e$ in a hyperbolic tetrahedron $\sigma$, and $l^+$ is defined in Definition \ref{def-l+}.
\end{definition}

\begin{rmk}
If the triangulation $\CT$ is not closed, for each boundary edge, $2\pi$ in the definition of its (generalized) curvature needs to be replaced with $\pi$. However, the proofs of our main theorem have no essential difference in this setting.
\end{rmk}

The curvature $K:\CL(M,\CT)\to{\IR}^{E}$ and the extended curvature $\wt{K}:{\IR}^{E}\to{\IR}^{E}$ are the focus of our research below. Note each decoration $w\in \IR^V$ acts linearly on ${\IR}^E$ by
\[(w+x)(vv') = w(v) + w(v') + x(vv')\]
where $x\in {\IR}^E$, and $vv'$ is an edge in $E$ with two vertices $v$ and $v'$. Given a partially truncated and decorated hyperbolic tetrahedron, if $v$ is a hyper-ideal vertex, we define $w(v) = 0$. For any two decorated hyperbolic polyhedral metrics $l_1$ and $l_2$ on $\CL(M,\CT)$, define an equivalent relationship $\sim$
on $\CL(M,\CT)$, that is, $l_1 \sim l_2$ if there exists a change of decoration $w \in \IR^V$ such that
$l_1 = l_2 + w$. Hence the quotient space  $\CL(M,\CT)/\sim$ is the space $\CL(M,\CT)$ modules the equivalent relationship $\sim$.

\subsection{Proof of the global rigidity}

Then, we can have the following result:

\begin{thm}\label{global rigidity-all}
For any triangulated compact pseudo 3-manifold $(M,\CT)$, a decorated hyperbolic polyhedral metric on $(M,\CT)$ is determined up to isometry and change of decorations by its Ricci curvature, that is, the curvature map
$K:\CL(M,\CT)/\sim \to {\IR}^{E}$ is injective.
\end{thm}

\begin{proof}

For each $l\in \CL(M,\CT)$ and each tetrahedron $ \sigma \in \CT$, let $l_{\sigma} \in \CL \subset \CR$ be the edge length vector of
$\sigma$ in the decorated hyperbolic polyhedral metric $l$.  Now we can define
\[ cov(l)=\sum_{ \sigma \in \CT} cov(l_{\sigma})  .\]

Based on Theorem \ref{cov-po}, Corollary \ref{general-cov}, Remark \ref{cov-mixed} and results in \cite{Luo2018}
and \cite{Feng-Ge-Liu}, for any partially truncated and decorated hyperbolic tetrahedron $\sigma$, the Hessian matrix of $cov(l_\sigma)$ is positive definite at each point in $\CL / \sim$, this means the Hessian matrix of $cov$ is positive definite. Thus, $cov$ is locally strictly convex in $\CL(M,\CT) / \sim$. On the other hand, 
by directly calculation, $\nabla cov = (2\pi,\ldots, 2\pi) - K_l $.

Now suppose otherwise that there exist $l_1, l_2 \in \CL({M,\CT})/\sim $ so that
$K_{l_1}=K_{l_2}$ but $l_1 \neq l_2$. Connecting $l_1$ and $l_2$ in ${\IR}^E$ by the line segment
$tl_1 + (1-t)l_2$, $t\in [0,1]$, and consider the function $m(t)=cov (tl_1 + (1-t)l_2) $, $ t\in [0,1] $.

By the construction, $m(t)$ is a $C^1$-smooth convex function so that $ m'(t)=\nabla cov \cdot (l_2-l_1) $. Since $\nabla cov(l_i) = (2\pi,\ldots,2\pi) - K_{l_i}$ and $K_{l_1} = K_{l_2}$.
Thus $m'(0)=m'(1)$. Because of the convexity of $m(t)$, $m(t)$ must be a linear function in $t$.
Since $l_1$ and $l_2$ are not equivalent, $cov$ is strictly convex near $l_1$ and $l_2$ along the line connecting $l_1$ and $l_2$.
Thus $m(t) $ is strictly convex in $t$ near 0 and 1. This contradicts that $m(t)$ is a linear function.

Therefore the curvature map $K:\CL(M,\CT)/\sim \rightarrow {\IR}^{E}$ is injective, that is,
for any $l_1 ,\; l_2 \in \CL({M,\CT}) $, if $K_{l_1}=K_{l_2}$, then there exists a change of decoration $w\in \IR^V$ such that
$ l_1 = l_2 + w $.
Hence, a decorated hyperbolic polyhedral metric on $(M,\CT)$ is determined up to isometry and change of decorations by its curvature.
\end{proof}

\begin{cor}\label{global-rigidity}
For any triangulated compact pseudo 3-manifold $(M,\CT)$, $k=2,3$, a decorated $k$-$(4-k)$ type hyperbolic polyhedral metric on $(M,\CT)$ is determined up to isometry and change of decorations by its curvature.
\end{cor}

Previously, Luo-Yang \cite{Luo2018} obtained the above rigidity for $k=0$ or $4$, and \cite{Feng-Ge-Liu} obtained the above rigidity for $k=1$. Comparing Theorem \ref{global rigidity-all} and Corollary \ref{global-rigidity}, we see the rigidity is valid not only for a space glued solely by one of the five type tetrahedra, but also for a space glued by maybe different type tetrahedra. Hence Corollary \ref{global-rigidity} parallels the rigidity obtained in \cite{Luo2018} and \cite{Feng-Ge-Liu}, while Theorem \ref{global rigidity-all} extends their rigidity to the most general setting. In addition, there had been many very important works on rigidity of hyperbolic cone metrics on 3-manifolds, like the works of Hodgson-Kerckhoff \cite{Hodgson1998}, Weiss \cite{Weiss2013}, Mazzeo-Montcouquiol \cite{Mazzeo2011} and Fillastre-Izmestiev \cite{Fillastre2009}\cite{Fillastre2011}\cite{Izmestiev} and others. As had been stated in \cite{Luo2018}, the difference between their work and ours, is that we consider the case where the singularity consists of complete geodesics from cusp to cusp or geodesics orthogonal to the totally geodesic boundary with possible cone singularities.

\subsection{The $H$-function}
\label{section-H-function}

\begin{definition} The functional $H:\IR^E\to\IR$ is defined as
 $$H(l) = cov(l) - 2\pi\sum_{i\in E} l_i,\quad l\in \IR^E.$$
\end{definition}

By Theorem~\ref{cov-po} and Corollary~\ref{general-cov}, we have the following result.
\begin{prop}\label{prop:stc} The functional $H$ is a $C^1$-smooth convex function on $\IR^E,$ which is smooth and locally strictly convex on $\CL(M,\CT) / \sim.$
\end{prop}

For a tetrahedron $\sigma\in\CT$ with the generalized decorated hyperbolic polyhedral metric given by $l,$  for any $e\prec\sigma$,
$$\frac{\partial cov_{\sigma}(l)}{\partial (l(e))}=\alpha_e(l).$$
Hence for any $j\in E,$
\begin{eqnarray}
\frac{\partial H}{\partial l_j}
&=& \frac{\partial cov(l)}{\partial l_j} - 2\pi
=\frac{\partial}{\partial l_j}\left(\sum_{j\prec\sigma} cov_{\sigma}(l)\right)-2\pi\label{eq:deriv}\\
&=& \sum_{j\prec\sigma}\alpha_{\sigma, e}(l) - 2\pi\nonumber\\
&=& -\wt{K}_j.\nonumber
\end{eqnarray}

\section{The combinatorial Ricci flow}
\label{section-CRF}
\subsection{Extended combinatoral Ricci flow}
\label{section-Extend-RCF}
Now, continuing with the setting of Chapter 5, we introduce the extended combinatorial Ricci flow in $\CL(M, \CT)$. It can be written as follows
\begin{equation}
\begin{cases}\label{eq:luoflow}
\frac{d l_i(t)}{d t} = \wt{K}_i(l(t)) ,\quad \forall i\in E, t\geq 0,\\
l(0) = l_0,
\end{cases}
\end{equation}
where $l_0\in \CL(M, \CT)$ and $l(t)\in \IR^{E(\CT)},\forall t>0.$ Since $\CL(M, \CT)$ is an open subset in $\IR^E$ and $K(l)$ is smooth on $\CL(M, \CT)$, the Picard theorem in the standard ODE theory yields the following.
\begin{thm}
For a triangulated compact pseudo 3-manifold $(M, \mathcal{T}),$ for any initial data $l_0\in\CL(M, \CT),$ the solution $\{l(t)|t\in [0, T )\}\subset \CL(M, \CT)$ to the extended combinatorial Ricci flow \eqref{eq:luoflow} exists and is unique on the maximal existence interval $[0, T )$ with $0 < T\leq \infty.$
\end{thm}

Next, using the same proof argument of Theorem 1.3 in \cite{Feng2022-1} we prove the long-time existence and uniqueness of the extended combinatorial Ricci flow. 
\begin{thm}\label{UNI}
For a triangulated compact pseudo 3-manifold $(M, \mathcal{T})$ and any initial data $l_0\in\IR^E,$ there exists a unique solution $\{l(t)|t\in [0,\infty)\}\subset \IR^E$ to the extended combinatorial Ricci flow \eqref{eq:luoflow}.
\end{thm}

By (\ref{eq:deriv}), the extended combinatorial Ricci flow (\ref{eq:luoflow}) is the negative gradient flow of the functional $H,$ and $l\in \CL(M,\CT)$ has zero Ricci curvature if and only if $l$ is a critical point of the functional $H.$

\begin{definition} Let $(M,\CT)$ be a triangulated compact pseudo 3-manifold. We say that the solution of the extended combinatorial Ricci flow is \emph{compatible with the change of decorations} (i.e. the action $\IR^V$ on $\IR^E$) if
for any initial data $l_0\in \IR^E$ with $l_0=w+l_0^\top,$ where $w\in \IR^V$ and $l_0^\top$ is the projection of $l_0$ to $\mathbb{R}^{E}/ \widehat{\mathbb{R}^{V}},$ then the solution $l(t)$ of the extended combinatorial Ricci flow with the initial data $l_0$ satisfies
$$l(t)=w+l^\top(t), \ \forall t\in[0,\infty),$$ where $l^\top(t)\in \mathbb{R}^{E}/ \widehat{\mathbb{R}^{V}}$ is the solution of  the extended combinatorial Ricci flow with the initial data $l_0^\top.$  Similarly, one defines the combinatorial Ricci flow is compatible with the change of decorations.
\end{definition}

Yang \cite{Y} proved that the combinatorial Ricci flow is compatible with the change of decorations for cusped 3-manifolds.
\begin{lemma}[Proposition~1.33 in \cite{Y}]\label{invariance of sum of edge length}
Let $(M, \mathcal{T})$ be a cusped $3$-manifold with an ideal triangulation.
Then for any $p\in V,$ $\sum_{i\in E, p\prec i}l_i$ is invariant along the combinatorial Ricci flow. In particular, the combinatorial Ricci flow is compatible with the change of decorations.
\end{lemma}

This is mainly because the sum of three inner angles in one decorated vertex triangle is $\pi$, for any ideal vertex. Thus, in the similar argument, we know:
\begin{prop}[Proposition~3.6 in \cite{Feng2022-1}]
\label{prop:inv}
Let $(M,\CT)$ be a pseudo 3-manifold with a closed triangulation which has a zero-curvature generalized decorated metric. Then the extended combinatorial Ricci flow is compatible with the change of decorations.
\end{prop}

\subsection{Fundamental Theorem of Combinatorial Ricci Flow}

\begin{prop}\label{H}
The functional $H$ is non-increasing along the extended combinatorial Ricci flow (\ref{eq:luoflow}), i.e. for any solution $l(t)$ to the flow (\ref{eq:luoflow}),
\[\frac{\partial H(l(t))}{\partial t} \le 0.\]
\end{prop}
\begin{proof}
By direct calculation,
\[\frac{\partial H(l(t))}{\partial t} = -|\wt{K}(l(t))|^2 \le 0.\]
\end{proof}
\begin{prop}\label{prop:zero}
If a solution $l(t)$ of the extended combinatorial Ricci flow \eqref{eq:luoflow} converges to some $\widetilde l \in \mathbb{R}^E$ as $t \to +\infty$. Then $\widetilde K(\widetilde l) = 0.$ In particular, if the solution $l(t)$ converges to $\bar l \in \CL(M, \CT)$. Then $K(\bar l) = 0$.
\end{prop}
\begin{proof}
 By Proposition \ref{H}, $H(l(t))$ is non-increasing. Moreover,
$\{H(l(t)) : t \ge 0\}$ is bounded, since $H$ is continuous on $\CL$ and $l(t) \to \widetilde  l$, $t \to \infty$.
Hence the following limit exists and is finite,
\[\lim_{t \to \infty} H(l(t)) = C.\]
Consider the sequence $\{H(l(n))\}^\infty_{n=1}$. By the mean value theorem, for any
$n \ge 1$ there exists $t_n \in (n, n + 1)$ such that
\begin{equation}
H(l(n + 1)) - H(l(n)) =\left.\frac{d}{dt} \right|_{t = t_n} H(l(t)) = - |\wt{K}(l(t_n))|^2.
\end{equation}
Note that $\lim_{n \to \infty} H(l(n + 1)) - H(l(n)) = 0$. We get
\[\lim_{n \to \infty} \wt{K}_i(l(t_n)) = 0, \quad \forall i\in E.\]
Since $l(t_n) \to \widetilde l$ as $n \to \infty$, the continuity of $\wt{K}_i$ yields that $\wt{K}_i(\widetilde  l) = 0$ for any $i\in E$.

Particularly, if $\widetilde l = \bar l \in \CL(M, \CT)$, then $\wt{K}_i(\bar l) = K_i(\bar l) = 0$.
\end{proof}


\begin{thm}\label{thm:l1}
Let $(M,\CT)$ be a compact pseudo 3-manifold with the closed triangulation which has a zero-curvature decorated metric $l \in \CL(M, \CT)$. Then $ H(l)$ is proper on $\mathbb{R}^{E} / \sim$, i.e.
 $$\lim_{l\in\mathbb{R}^{E}/\sim,l \to \infty}  H(l) = +\infty.$$
\end{thm}
\begin{proof}

We consider the functional $H:\IR^E\to \IR.$ By Proposition~\ref{prop:stc}, it is $C^1$-smooth and convex on $\IR^E,$ and is smooth and strictly convex on $\mathcal{L}(M, \mathcal{T})/\sim.$ By \eqref{eq:deriv}, the metrics in $\IR^E$ with zero-curvature correspond to the critical points of the functional $H.$ Note that any critical point of a convex function on a convex domain is a minimizer. So that the zero-curvature metric $l$ is a minimizer of $H.$ Moreover, $l\in \mathcal{L}(M, \mathcal{T})$ and $H$ is strictly convex on $\mathcal{L}(M, \mathcal{T})/\sim.$ This is a unique minimizer, also a unique critical point, of $H$ on $\IR^E /\sim.$  By the strict convexity on $\mathcal{L}(M, \mathcal{T})/\sim,$ we have $$\lim_{l\in\mathbb{R}^{E}/ \sim,\; l\to\infty} H(l)=+\infty.$$
\end{proof}

\begin{thm}\label{thm:main1}
For a compact pseudo 3-manifold with a closed triangulation $(M, \mathcal{T}),$ let $\{l(t)\}_{t \ge 0}$ be a solution to the extended combinatorial Ricci flow \eqref{eq:luoflow}.
\begin{enumerate}
\item If there is no zero-curvature generalized decorated metrics, then $l(t)$ diverges to infinity in subsequence, i.e. there exists a subsequence $t_n\to\infty$ such that $|l(t_n)|\to \infty, n\to\infty.$
\item There is a zero-curvature decorated metric if and only if $l(t)$ converges to a zero-curvature decorated metric for any initial data. In this case, the convergence is exponentially fast.
\end{enumerate}
\end{thm}

\begin{proof}[Proof of Theorem~\ref{thm:main1}]
Let $\phi(t) = H(l(t)).$

We prove the first assertion. Suppose that it is not true, then there exists a constant $C$ such that
$|l(t)|\leq C,\ t\in [0,\infty).$ 
Then $\phi(t)$ is bounded on $t\in[0,\infty).$ Since $\phi(t)$ is non-increasing,  the following limit exists
$$\lim_{t\to\infty} \phi(t)=C\in \IR.$$ 
By the mean value theorem, there are $\xi_n\in (n,n+1)$ such that $$\phi'(\xi_n) = \phi(n + 1) - \phi(n) \to 0,\quad n\to\infty$$
Since $\phi'(t)= - |\widetilde K(l(t))|^2,$ $$\wt K(l(\xi_n))\to 0,\quad n\to\infty.$$
There exists a subsequence of $\{\xi_n\}_{n=1}^{\infty}$, denoted by $\{\xi_{n_k}\}_{k=1}^{\infty}$, such that $l(\xi_{n_k})\to l_\infty\in \IR^E.$ This implies that $$\wt{K}(l_\infty)=0.$$ This yields a contradiction and proves the first assertion.

~

Now we prove the second assertion. By Proposition~\ref{prop:zero}, we only need to prove that if there is a zero-curvature decorated metric, then $l(t)$ converges to a zero-curvature decorated metric exponentially fast. By Proposition~\ref{prop:inv}, it suffices to prove the result for the initial data $l_0\in \mathbb{R}^{E}/ \sim.$ Hence the solution $l(t)$ of the extended combinatorial Ricci flow  satisfies $l(t)\in \mathbb{R}^{E}/ \sim $ for all $t\in[0,\infty).$

By Theorem~\ref{global rigidity-all}, there is a unique zero-curvature decorated metric $\hat{l}$ in
$\mathbb{R}^{E}/ \sim.$
By Theorem~\ref{thm:l1}, the functional $H|_{\mathbb{R}^{E}/ \sim }$ is proper on $\mathbb{R}^{E}/ \sim $ and $H$ is bounded from below.


For $\phi(t) = H(l(t)),$ since $\phi(t)$ is non-increasing, the following limit exists
$$\lim_{t\to\infty} \phi(t)=C\in \IR.$$ Combining the properness of $\wt{H}|_{\mathbb{R}^{E}/ \sim },$ $\{l(t):t\in [0,\infty)\}$ is contained in a compact subset of $\mathbb{R}^{E}/ \sim .$
By the mean value theorem, there is a sequence $\xi_n\in (n,n+1)$ such that
$$-|\wt{K}(l(\xi_n))|^2=\phi'(\xi_n) = \phi(n + 1) - \phi(n) \to 0,\quad n\to\infty.$$

Passing to a subsequence, still denoted by $\xi_n,$
$$l(\xi_n)\to l_\infty,\quad n\to\infty.$$
By the continuity of $\wt{K},$  we have $\wt{K}(l_\infty)=0.$
By Theorem~\ref{global rigidity-all}, we have $l_\infty=\hat{l}.$

Hence for any neighbourhood $U$ of $\hat{l}$ in $\mathbb{R}^{E}/ \sim,$ for sufficiently large $\xi_n,$
$l(\xi_n)\in U.$
By Proposition~\ref{prop:stc},
$Hess (\wt{H}|_{\mathbb{R}^{E}/ \sim})_{\hat{l}}$ is positive definite, which implies the critical point $\hat{l}$ is a local attractor of the extended combinatorial Ricci flow restricted on $\mathbb{R}^{E}/ \sim.$ By Lyapunov Stability Theorem, \cite[Ch.5]{Pon}, $l(t)$ converges to $\hat{l}$ exponentially fast for any initial data. This proves the result.
\end{proof}

Finally, suppose $M$ is a oriented compact 3-manifold with boundary, no component of which is a $2$-sphere, and $\CT$ is an (topological) ideal triangulation of $M-\partial_t$, where $\partial_t$ is the toral boundary components of $M$. We remark that $M$ always admits such a triangulation, see Bing \cite{Bing}, Moise \cite{Moise} and Jaco-Rubinstein \cite{Jaco} (or see \cite{GeJZ-decom,GeJZ-angle}). Hence we obtain a compact pseudo 3-manifold $M$ with a closed triangulation $\CT$. Consider the decorated hyperbolic polyhedral metrics on $(M,\CT)$, each metric $l\in\CL(M, \CT)$ produces a hyperbolic cone metric on $M$.
By Definition \ref{dhpm}, if there is a metric $l\in\CL(M,\CT)$ with zero Ricci curvature at each edge, then we obtain a complete hyperbolic metric with totally geodesic boundary on $M$, and furthermore, $\CT$ is a geometric ideal triangulation of the hyperbolic manifold $M-\partial_t$. Hence the fundamental theorem of CRF, i.e. Theorem \ref{Fundamental-Thm} follows naturally from Theorem \ref{thm:main1}.


\subsection{Some comments to Corollary \ref{thm-converge-imply-topo}}

Suppose $M$ is an oriented compact 3-manifold with boundary, no component of which is a 2-sphere, $\CT$ is an ideal triangulation of $M$. If $(M,\CT)$ admits an angle structure (see \cite{Futer2011,GeJZ-angle} for a definition), then by Lackenby's work \cite{Lackenby2000,Lackenby-AGT}, the topology of the manifold $M$ is restricted, that is, $M$ is irreducible, atoroidal and not Seifert fibred. By Thurston's hyperbolization theorem, $M$ is hyperbolic. On the contrary, given a hyperbolic 3-manifold $M$, it always admits an ideal triangulation. However, whether it admits an angle structure or not is still an open question. Of course, if Thurston's geometric ideal triangulation conjecture holds, then there is an angle structure. Equivalently, if the solution $l(t)$ to the extended CRF converges, then the limit metric $l(+\infty)$ determines an angle structure, so that any connected closed 2-normal surface in $\CT$ with non-negative Euler characteristic is normally parallel to a toral boundary components of $M$.

$0$-efficient and $1$-efficient triangulations of 3-manifolds were introduced and extensively studied in Jaco-Rubinstein \cite{Jaco}, Kang-Rubinstein \cite{Kang} and  Garoufalidis-Hodgson-Rubinstein-Segerman \cite{Garoufalidis}, which have applications to finiteness theorems, knot theory, Dehn fillings, decision problems, algorithms, computational complexity and Heegaard splittings. By definition, an ideal triangulation of an orientable 3-manifold $M$ is 0-efficient if there are no embedded normal 2-spheres or one-sided projective planes. In addition, $\CT$ is 1-efficient if it is 0-efficient, the only embedded normal tori are vertex-linking and there are no embedded one-sided normal Klein bottles. $\CT$ is strongly 1-efficient if there are no immersed normal 2-spheres, projective planes or Klein bottles and the only immersed normal tori are coverings of the vertex-linking tori. The proof of Theorem 1.5 in \cite{Garoufalidis}, Theorem 2.5 and Theorem 2.6 in \cite{Kang} can be applied almost word for word to our situation. So we get our Corollary \ref{thm-converge-imply-topo}.

In the pioneering work of Costantino-Frigerio-Martelli-Petronio \cite{CFMP}, they first use the edge valence to obtain a hyperbolic structure on $(M,\CT)$. Here the edge valence of an edge $e\in E$ is the number of tetrahedra surrounding $e$. They showed that if all valencies are no less than 6, then $(M,\CT)$ admits a hyperbolic structure. Moreover, they conjectured that the original triangulation $\CT$ is also geometric. The authors and Hua \cite{Feng2022} (or see recent progress in \cite{Feng-Ge-Meng}) proved their conjecture for any triangulation with valencies no less than $10$, by showing that the extended CRF converges. Thus a direct way to control the extended CRF, and hence then prove Thurston's geometric triangulation conjecture, is to find a topological ideal triangulation with edge valence no less than $10$. On the other hand, maybe we can use the CRF with surgery to find suitable triangulations. In dimension 2, combinatorial Yamabe flows with surgery first appeared in \cite{Gu-Luo-I,Gu-Luo-II} with the help of Teichm\"{u}ller theory. In this case, the surgery is performed by replacing one diagonal of a quadrilateral (made up of two triangles pasted along a common edge) with another, and the time of the operation is exactly when the four points are in a circle, so as to ensure that the triangulation is Delaunay. These methods should not be easy to generalize to three dimension, so we don't know how to perform surgery on three dimensional CRF at this time. It was shown in \cite{Jaco} that any ideal triangulation $\CT$ can be modified to an $0$-efficient ideal triangulation. This modification can be achieved by the so called crushing operation. It is conjectured that Thurston's geometric ideal triangulation (or the $0$-efficient ideal triangulation, respectively) should correspond to a certain stability, both from the perspective of the stability of the combinatorial structure and from the perspective of a certain energy minimization. The CRF with surgery may be suitable to find Thurston's geometric decomposition (or achieve this modification, respectively) since it tends to minimize the $H$-function.

\noindent Ke Feng, kefeng@uestc.edu.cn\\[2pt]
\emph{School of Mathematical Sciences, University of Electronic Science and Technology of China, Sichuan 611731, P. R. China}\\[2pt]

\noindent Huabin Ge, hbge@ruc.edu.cn\\[2pt]
\emph{School of Mathematics, Renmin University of China, Beijing 100872, P. R. China}\\[2pt]

\end{document}